\documentclass[12pt]{amsart}

\usepackage{amssymb,amsthm,amsmath,amscd}

\usepackage{mathtools}
\usepackage{soul}
\usepackage{booktabs}
\RequirePackage[dvipsnames,usenames]{color}

\usepackage[hidelinks]{hyperref}

\input{kmacros3.sty}

\usepackage{appendix}
\usepackage{tikz}
\usepackage{tikz-cd}
\usetikzlibrary{cd}
\usepackage{graphicx}
\usepackage[all,cmtip]{xy}

\usepackage{mabliautoref}
\numberwithin{equation}{theorem}

\usepackage{bm}
\usepackage{pifont}
\usepackage{upgreek}

\usepackage{eucal}
\usepackage{ulem}
\usepackage{stmaryrd}

\usepackage{fullpage}

\usepackage{setspace}
\usepackage{enumerate}
\usepackage{calc}

\usepackage{bbding}

\usepackage{verbatim}
\usepackage{alltt}
\usepackage[capitalise]{cleveref}
\usepackage[dvipsnames]{xcolor}
\definecolor{chianti}{rgb}{0.6,0,0}
\definecolor{meretale}{rgb}{0,0,.6}
\definecolor{leaf}{rgb}{0,.35,0}

\theoremstyle{plain}
\newtheorem{maintheorem}{Main Theorem}

\usepackage[most]{tcolorbox}

\newtcolorbox{commentbox}{
    enhanced,
    breakable,
    colback=white,
    colframe=OliveGreen,
    colbacktitle=OliveGreen,
    coltitle=white,
    fonttitle=\bfseries,
    title=Comment,
}

\newtcolorbox{todobox}{
    enhanced,
    breakable,
    colback=white,
    colframe=Mahogany,
    colbacktitle=Mahogany,
    coltitle=white,
    fonttitle=\bfseries,
    title=To Do,
}

\DeclareMathOperator{\vol}{vol}

\newcommand{\fm}{\mathfrak{m}}

\newcommand{\QQ}{\mathbb{Q}}
\newcommand{\CC}{\mathbb{C}}

\newcommand{\NN}{\mathbb{N}}
\newcommand{\fa}{\mathfrak{a}}

\newcommand{\FF}{\mathbb{F}}

\DeclareMathOperator{\edim}{edim}

\DeclareMathOperator{\Span}{Span}

\DeclareMathOperator{\lct}{lct}
\DeclareMathOperator{\FPT}{fpt}

\begin{document}

\title{Uniform Positivity of $F$-Signature and Effective Ideal Containments}

\author{Hongyu Luo}
\address{Department of Mathematics, University of Alabama, Tuscaloosa, AL 35487 USA}
\email{hluo13@crimson.ua.edu}

\author{Thomas Polstra}
\address{Department of Mathematics, University of Alabama, Tuscaloosa, AL 35487 USA}
\email{tmpolstra@ua.edu}
\thanks{Polstra was supported by a grant from the National Science Foundation, DMS 2502317.}

\begin{abstract}
We establish inequalities between numerical invariants and ideal containments in local rings containing a field. We prove the positivity of $F$-signature for klt-type algebras essentially of finite type over $\mathbb{C}$ under prime characteristic reduction and resolve the Uniform Symbolic Topology Problem effectively for both $F$-finite strongly $F$-regular rings of prime characteristic and klt type algebras over $\mathbb{C}$. Consequently, letting $s(R)$ denote the $F$-signature and $d$ the Krull dimension, if $d\geq 2$ then $I^{\left(\left\lceil\frac{d-1}{s(R)}\right\rceil n\right)} \subseteq \overline{I^{n}}$ and $I^{\left(d\left\lceil\frac{d-1}{s(R)}\right\rceil n\right)} \subseteq I^n$ for every ideal $I \subseteq R$ and $n \in \mathbb{N}$ in these settings.
\end{abstract}

\maketitle

\section{Introduction}

Throughout this article, $R$ is a commutative noetherian ring and $(R,\fm,k)$ denotes a local commutative noetherian ring with unique maximal ideal $\fm$ and residue field $k$.

\subsection{$F$-signature and $\alpha_F$-invariant}

In this subsection, let $(R,\fm,k)$ be a local domain of prime characteristic $p>0$, $F:R\to R$ be the Frobenius endomorphism, and for each $e\geq 1$ let $F^e_*R$ be the $R$-module obtained by restriction of scalars along the $e$th iterate of the Frobenius endomorphism. Assume further that $R$ is $F$-finite so that $F^e_*R$ is a finitely generated $R$-module for every $e\in\NN$.

We denote by $a_e(R)$ the maximal rank of a free $R$-summand of $F^e_*R$. The local $F$-finite domain $(R,\fm,k)$ is regular if and only if $a_e(R) = \rank_R(F^e_*R)$ for some, equivalently for all $e\geq 1$, \cite{Kunz69}. The ring $R$ is $F$-pure if $a_e(R)\geq 1$ for some, equivalently every, positive integer $e$. The local ring $R$ is strongly $F$-regular if for each $0\not=x\in R$ there exists $e\in\NN$ so that $F^e_*x$ generates a free $R$-summand of $F^e_*R$.

For each $e\geq 1$ let $I_e(\fm) = \{x\in R\mid R\xrightarrow{\cdot F^e_*x}F^e_*R\mbox{ does not split}\}$ denote the $e$th splitting ideal of $R$. Then $x\in R\setminus I_e(\fm)$ if and only if $F^e_*x$ generates a free $R$-summand of $F^e_*R$. The $F$-pure threshold of $x$ asymptotically measures the free summands of $F^e_*R$ being contributed by powers of $x$. Specifically, if $\nu_e(x) = \max\{t\in\NN\mid x^t\not\in I_e(\fm)\}$, the $F$-pure threshold of $x$ is characterized as the limit $\fpt_R(x):= \lim_{e\to \infty} \frac{\nu_e(x)}{p^e}$.

The $F$-signature and $\alpha_F$-invariant of a local $F$-finite domain of prime characteristic are complementary numerical invariants detecting strong $F$-regularity. The $F$-signature of $R$, denoted by $s(R)$, is the asymptotic measurement of the largest rank of a free summand of $F^e_*R$ and defined as
\[
s(R) = \lim_{e\to \infty}\frac{a_e(R)}{\rank_R(F^e_*R)},
\]
a well-defined limit by \cite{Tucker12}. Recently introduced and developed in \cite{Pande23,TY26}, the $\alpha_F$-invariant of an $F$-finite and $F$-pure domain $(R,\fm,k)$ is defined as
\[
\alpha_F(R) = \inf\{\fpt_R(x)\overline{\ord}_{\fm}(x)\mid 0\not=x\in \fm\} 
\]
where $\overline{\ord}_{\fm}(x)$ is the normalized order of $x$.

Both $F$-signature and the $\alpha_F$-invariant take values between $0$ and $1$, and are positive if and only if $R$ is strongly $F$-regular, \cite[Main Result]{AL03} and \cite[Proposition~3.9 and Remark~3.11]{TY26}. Positivity of $\alpha_F(R)$ being equivalent to strong $F$-regularity is derived from the inequalities of \cite[Proposition~3.9]{TY26}. Alternatively, we prove quite generally that $\alpha_F(R)\geq s(R)$ in \cref{cor: F-sig lower bound for alpha F invariant}, a result which allows positivity of $\alpha_F(R)$ being equivalent to strong $F$-regularity to be more directly derived.

Our main contributions to $F$-singularity theory are lower bounds for the $\alpha_F$-invariant and the $F$-signature of a local ring of prime characteristic that depend only upon the dimension of $R$, the Hilbert-Samuel multiplicity of $R$, a constant derived from the Izumi-Rees Theorem, and the $F$-pure threshold of a suitably chosen Jacobian element.

\begin{maintheorem}
    \label{main theorem: F-sig lowerbound via fpt of single element}
    Let $(R,\fm,k)$ be a local strongly $F$-regular $F$-finite domain of prime characteristic $p>0$ and Krull dimension $d$. Assume that $R$ is not regular. Let $(A,\fm_A,k)\subseteq (R,\fm,k)$ be a regular local subring of $R$ so that $R$ is finite and generically \'{e}tale over $A$. Suppose that $z_1,\ldots,z_c\in R$ generate $R$ as an $A$-algebra, and $R\cong A[z_1,\ldots,z_c]/P$. Let $f_1,\ldots,f_c\in P\subseteq A[z_1,\ldots,z_c]$ be local generators of $P$, $D = \det(\frac{\partial f_i}{\partial z_j})_{1\leq i,j\leq c}$, and $\delta$ the (necessarily nonzero) image of $D$ in $R$. If $E$ is a constant so that for each $0\not=x\in \fm$, 
        \[
        e_{\frac{(\fm_A,x)R}{xR}}(R/xR)\leq E\overline{\ord}_{\fm}(x),
        \]
        where $\overline{\ord}_{\fm}(x)$ is the normalized order of $x$, then
        \[
        \alpha_F(R)\geq \frac{\fpt_R(\delta)}{2E}\quad\mbox{and}\quad s(R) \geq \frac{e(R)}{d!}\left(\frac{\fpt_R(\delta)}{2E}\right)^d.
        \]
\end{maintheorem}

\begin{remark}
    The existence of the constant $E$ in Main Theorem~\ref{main theorem: F-sig lowerbound via fpt of single element} follows from the Izumi-Rees theorem, which establishes the linear comparability of divisorial valuations centered at the maximal ideal of an analytically irreducible local domain \cite{Rees1989}, together with Rees's degree function formula \cite{Rees1961}. If $\overline{\fm_A R}=\fm$ and $\Gr_{\fm}(R)$ is a normal domain, then we may take $E=e(R)$. These conditions hold, for example, when $R$ is the completion of a standard graded normal domain at its homogeneous maximal ideal and $A$ is the completion of a Noether normalization generated by degree-one elements.
\end{remark}

\subsection{Limit of $F$-signature}

We now turn our attention to the behavior of $F$-signature in mod $p$ reductions of a normal local domain $(R_{\CC},\fm,k)$ essentially of finite type over $\CC$. Let $X_{\CC}=\Spec(R_{\CC})$ and $K_{X_{\CC}}$ be a canonical divisor. We say that $R_{\CC}$ is of klt type if there is an effective $\QQ$-divisor $\Delta$ on $X_{\CC}$ such that $K_{X_{\CC}}+\Delta$ is $\QQ$-Cartier and the pair $(X_{\CC},\Delta)$ is Kawamata log terminal (klt). The definition of klt type does not require $K_{X_{\CC}}$ itself to be $\QQ$-Cartier and this formulation is compatible with the treatment of singularities on non-$\mathbb{Q}$-Gorenstein varieties in \cite{deFernexHacon2009}.

To motivate mod $p$ reductions in the introduction, assume for the sake of simplicity $R_{\CC}$ is realized as follows. Let $\ZZ\subseteq R = \ZZ[x_1,\ldots,x_n]/P$, where $P$ is a prime ideal contained in $(x_1,\ldots,x_n)\ZZ[x_1,\ldots,x_n]$, and $R_{\CC} = \left(\CC[x_1,\ldots,x_n]/P\right)_{(x_1,\ldots,x_n)}$. For each prime number $p>0$ let $R_p = \left(\FF_p[x_1,\ldots,x_n]/P\right)_{(x_1,\ldots,x_n)}$. Then the correspondence of singularities of the complex minimal model program and prime characteristic singularities provides that $R_{\CC}$ of klt type implies $R_p$ is strongly $F$-regular for all $p\gg 0$, \cite{HaraWatanabe2002, Takagi2004}.

If $R_{\CC}$ is of klt type, then results of \cite{MTW05, CEMS18,TY26,Takagi2004} and standard spreading out arguments give a uniform positive lower bound for the quantities in Main Theorem ~\ref{main theorem: F-sig lowerbound via fpt of single element} to be bounded from below. Consequently, we answer affirmatively \cite[Question~5.9]{CarvajalRojasSchwedeTucker2018}, more generally \cite[Conjecture A]{TY26} in \cref{cor: positive lim F-signature}, and several notable consequences indicated in \cite[Section~5]{CarvajalRojasSchwedeTucker2018}.

\begin{maintheorem}[Application of \cref{cor: positive lim F-signature}]
    \label{main thm: positive F-signature for klt type singularity}
    Let $\ZZ\subseteq R = \ZZ[x_1,\ldots,x_n]/P$ be a finitely generated normal domain over $\ZZ$, $P\subseteq (x_1,\ldots,x_n)$, and $R_{\CC}:= (\CC[x_1,\ldots,x_n]/P)_{(x_1,\ldots,x_n)}$. For each prime number $p$ let $R_p = (\FF_p[x_1,\ldots,x_n]/P)_{(x_1,\ldots,x_n)}$. If $R_{\CC}$ is of klt type then
    \[
    \liminf_{p\to\infty}s(R_p)>0.
    \]

    More generally, if $(R,\fm,k)$ is a local normal domain essentially of finite type over $\CC$, of klt type, and $s(R)$ the $F$-signature of $R$ as defined in \cite[Definition~4.1]{TY26}. Then $s(R)>0$.
\end{maintheorem}

The normalized volume $\widehat{\vol}(x,X,\Delta)$ of a pair $(X,\Delta)$ at a closed point $x\in X$, introduced by Li in \cite{Li18}, is a numerical invariant motivated by the study of K-stability and the existence of K\"ahler-Einstein metrics on Fano varieties. As an application of our lower bounds for $F$-signature in terms of $F$-pure thresholds of Jacobian elements described in Main Theorem~\ref{main theorem: F-sig lowerbound via fpt of single element}, when paired with valuation estimates of \cite{Zhuang24}, we establish a lower bound for the characteristic $0$ $F$-signature of a standard graded normal domain depending only upon the dimension, embedding dimension, and the normalized volume of a klt pair.

\begin{maintheorem}
    \label{mainthm: F-sig lower bound via normalized volume}
    Let $R$ be the localization of a standard graded normal domain over $\CC$ at its homogeneous maximal ideal. Let $d=\dim(R)$ and $d+c = \text{embdim}(R)$. Let $\Delta\geq 0$ be a boundary $\QQ$-divisor of $X = \Spec(R)$ and $x\in X$ the closed point. There exists constant $C(d,c)$ depending only upon $d$ and $c$ so that if $(X,\Delta)$ is klt then
    \[
    s(R) \geq C(d,c)\left(\widehat{\vol}(x,X,\Delta)\right)^d.
    \]
\end{maintheorem}


\subsection{Ideal Containments, Multiplicity Bounds, and Izumi-Type Inequalities}

In addition to the main results discussed thus far, we derive some effective ideal containment properties, multiplicity bounds, and Izumi-type results of local $F$-finite strongly $F$-regular domains and local normal domains essentially of finite type over $\CC$ of klt type. Notably, we recover, improve upon, and expand to the klt case a result of \cite{PolstraUSTP}: there exists a constant $C$ so that for every ideal $I\subseteq R$ and $n\in\NN$, $I^{(Cn)}\subseteq I^n$. Recall that if $I\subseteq R$ is an ideal, then the $n$th symbolic power of $I$ is $I^{(n)}:= I^nR_W\cap R$ where $R_W$ is the localization of $R$ at the complement of the union of the associated primes of $I$. The methods of \cite{PolstraUSTP} only show the existence of such a constant whereas here we can provide an explicit constant depending only upon dimension and the value of $F$-signature. 

\begin{maintheorem}
    \label{main theorem: Effective Uniform Bounds}
    Let $(R,\fm,k)$ be either a local $F$-finite and strongly $F$-regular domain of prime characteristic $p>0$ or a local normal domain essentially of finite type over $\CC$ and of klt type. Let $R$ have Krull dimension $d\geq 1$.
    \begin{enumerate}
        \item\label{main theorem: Effective Uniform Bounds USTP} Let $0\not=I\subseteq R$ be an ideal, $h$ the maximal height of an associated prime of $I$, and $n\in\NN$. Then
        \[
        I^{\left(\left\lceil\frac{h}{s(R)}\right\rceil n\right)}\subseteq \overline{I^n} \quad \mbox{and} \quad I^{\left(d\left\lceil\frac{h}{s(R)}\right\rceil n\right)}\subseteq I^n.
        \]
        \item\label{main theorem: Effective Uniform Bounds HS multiplicity bound} $e(R)\leq \frac{d!}{s(R)^{d-1}}$.
        \item\label{main theorem: Effective Uniform Bounds hypersurface multiplicity} If $I\subseteq R$ is an $\fm$-primary ideal, $0\not=x\in \fm$, then
        \[
        e_{\frac{(I,x)R}{xR}}(R/xR)\leq \frac{de_I(R)}{s(R)}\overline{\ord}_{\fm}(x).
        \]
        \item\label{main theorem: Effective Uniform Bounds Izumi-Rees} If $I$ is an $\fm$-primary ideal, $\nu_1,\nu_2$ are Rees valuations of $I$, and $0\not=x\in R$,
        \[
        \nu_1(x)\leq \frac{de_I(R)}{s(R)}\nu_2(x).
        \]
    \end{enumerate}
\end{maintheorem}

\begin{remark}
    The exponent of $F$-signature appearing in (\ref{main theorem: Effective Uniform Bounds HS multiplicity bound}) is optimal by $F$-signature computations of Singh, \cite{Singh2005}. To the best of the author's knowledge, the upper bound of the multiplicity of a principal quotient described in (\ref{main theorem: Effective Uniform Bounds hypersurface multiplicity}) and the linear comparison of Rees valuations of (\ref{main theorem: Effective Uniform Bounds Izumi-Rees}) is novel even in the case $R$ is regular, in which case $s(R) = 1$.
\end{remark}


\subsection*{AI Acknowledgments} The authors used ChatGPT (GPT-6 Astra) during the preparation of this article. Its main mathematical contribution was the development of \cref{thm: splitting ideals and valuation rings}. The authors prompted the model to analyze the methods of \cite[Lemma~13.19]{MaPolstra} and then to derive an explicit bound on the valuation of an element of a splitting ideal. In response, the model produced the statement and proof of \cref{thm: splitting ideals and valuation rings}.

Experts will be aware that Main Theorems~\ref{main theorem: F-sig lowerbound via fpt of single element} and~\ref{main thm: positive F-signature for klt type singularity} follow from \cref{lem: multivariable calculus lemma} and its prime characteristic corollary, \cref{cor: mapping to discriminant power}, via standard Frobenius techniques and reduction to positive characteristic arguments. Motivated by the prime characteristic techniques in \cite{PolstraSmirnov2020,DeStefaniSmirnov2024}, and historically rooted in methods of \cite{Samuel1951Singularites}, the authors independently established \cref{lem: multivariable calculus lemma} and \cref{cor: mapping to discriminant power} for hypersurfaces by proving the existence of a solution to a suitable power series equation. The authors developed this argument to supplement their earlier prime characteristic results in Main Theorem~\ref{main theorem: Effective Uniform Bounds} and without the intention of establishing results in characteristic $0$. To expedite the generalization process, the model assisted in extending the proof of \cref{lem: multivariable calculus lemma} to the general case by pointing out that the implicit function theorem readily generalized the power series argument provided in the hypersurface case.

The authors also used the model to help identify typographical errors and to suggest editorial revisions. The authors take responsibility for the mathematical content of the article.

\subsection{Overlap with Independent Work}
Ajit and Apsley~\cite{AjitApsley} independently established Main Theorem~\ref{main thm: positive F-signature for klt type singularity}. Their positivity argument is related to ours and \cite[Theorem~2.11]{AjitApsley} is comparable to our Main Theorem~\ref{main theorem: F-sig lowerbound via fpt of single element}. Their arguments rely upon clearing denominators of images of generalized Hasse-Schmidt differential denominators by powers of a distinguished element, ours uses a Jacobian determinant and a formal implicit function argument in a complete $F$-finite domain. Not obtained in this article, are methods that allow for the converse of Main Theorem~\ref{main thm: positive F-signature for klt type singularity}, whereas in \cite{AjitApsley} it is shown that uniform positivity of $F$-signature under reduction modulo $p$ characterizes singularities of klt type. Not present in \cite{AjitApsley}, our article develops additional quantitative relationships with local $\alpha_F$-invariants, multiplicity, and valuations. These additional developments provide several notable distinct developments. We give an effective solution to the Uniform Symbolic Topology Problem for $F$-finite strongly $F$-regular rings of prime characteristic and normal local domains essentially of finite type over $\CC$ of klt type. We also provide a uniform lower bound of characteristic $0$ $F$-signature that involves only the normalized volume, dimension, and embedding dimension of a klt pair in the standard graded setting.

\section*{Acknowledgments}
The authors thank Karl Schwede for informing us of the independent and overlapping work of \cite{AjitApsley}.

\section{Preliminaries}

\subsection{Multiplicity, Integral Closure, and Izumi-Rees}

Let $(R,\fm,k)$ be a noetherian local ring of Krull dimension $d$. If $I\subseteq R$ is an ideal, $I^n$ denotes the $n$th power of the ideal and set $I^n=R$ if $n\leq 0$. We denote by $\overline{I}$ the integral closure of an ideal $I\subseteq R$. If $I$ is $\fm$-primary, we denote its Hilbert-Samuel multiplicity by $e_I(R)$ and set $e(R)=e_{\fm}(R)$.

Throughout our discussion of Rees valuations, we assume that $R$ is an excellent analytically irreducible local domain. For a nonzero proper ideal $I\subseteq R$, let $\mathcal{R}_I$ denote the finite collection of Rees valuations of $I$. These are the divisorial valuations of $R$ corresponding to the exceptional components of the normalized blowup of $I$.

For $0\not=x\in R$, let $\ord_I(x)=\max\{n\geq 0\mid x\in I^n\}$ and $\overline{\ord}_I(x)=\lim_{n\to\infty}\ord_I(x^n)/n$. The Rees valuative criterion gives $\overline{\ord}_I(x)=\min_{\nu\in\mathcal{R}_I}\nu(x)/\nu(I)$ and, consequently, $x\in\overline{I^n}$ if and only if $\overline{\ord}_I(x)\geq n$; see \cite[Chapter~10]{HunekeSwanson2006}. Moreover, $\overline{\ord}_I(x^n)=n\overline{\ord}_I(x)$ and $\overline{\ord}_I(xy)\geq\overline{\ord}_I(x)+\overline{\ord}_I(y)$. For any discrete valuation $\nu$ centered on $\fm$, we also have $\nu(I)\overline{\ord}_{I}(x)\leq\nu(x)$.

Since $R$ is analytically unramified, for every ideal $I$ there exists an integer $b\geq 0$ such that $\overline{I^{n+b}}\subseteq I^n$ for all $n\geq 1$ by \cite[Theorem~9.1.2]{HunekeSwanson2006}. In particular, if $I$ is $\fm$-primary, then $\lim_{n\to\infty}d!\mathcal{L}_R(R/\overline{I^n})/n^d=e_I(R)$.

Let $I$ be $\fm$-primary and write $\mathcal{R}_I=\{\nu_1,\ldots,\nu_t\}$. By Rees's multiplicity formula \cite{Rees1961}, there exist positive integers $d_1,\ldots,d_t$, independent of $x$, such that for every $0\not=x\in\fm$,
\[
    e_{\frac{(I,x)R}{xR}}(R/xR)=\sum_{i=1}^t d_i\nu_i(x)
    \qquad\text{and}\qquad
    e_I(R)=\sum_{i=1}^t d_i\nu_i(I).
\]
See also \cite[Corollary~2.9]{Polstra2026} and \cite[Theorem~1.1]{CutkoskyMontano2025}.

Finally, the Izumi-Rees theorem asserts that any two divisorial valuations $\nu_1,\nu_2$ centered on $\fm$, there exists $C>0$ so that for any $0\not=x\in R$, $\nu_1(x)\leq C\nu_2(x)$; see \cite{Rees1989}. Applying this comparison to the finite collection $\mathcal{R}_I\cup\mathcal{R}_{\fm}$ and using the preceding formulas, there exists a constant $E=E(R,I)>0$ such that for every $0\not=x\in\fm$, 
\[
e_{\frac{(I,x)R}{xR}}(R/xR)\leq E\overline{\ord}_{\fm}(x).
\]

\subsection{$F$-singularities and Prime Characteristic Methods}

Throughout this subsection, $R$ is a noetherian ring of prime characteristic $p>0$. For $e\geq 1$, let $F^e:R\to R$ denote the $e$th iterate of the Frobenius endomorphism, and let $F^e_*M$ denote the $R$-module obtained from an $R$-module $M$ by restriction of scalars through $F^e$. The ring $R$ is $F$-finite if $F^e_*R$ is finitely generated for some, equivalently every, $e\geq 1$. An $F$-finite ring is $F$-pure if the Frobenius map $R\to F_*R$ is pure, equivalently split. An $F$-finite domain is strongly $F$-regular if, for every $0\not=x\in R$, there exists $e\geq 1$ such that the $R$-linear map $R\to F^e_*R$ sending $1$ to $F^e_*x$ splits.

For an $F$-finite ring $R$, an ideal $I\subseteq R$, and $e\geq 1$, the \emph{$e$th splitting ideal of $I$} is $I_e(I)=\{r\in R\mid \varphi(F^e_*r)\in I \text{ for every }\varphi\in\Hom_R(F^e_*R,R)\}$. Basic properties of splitting ideals are collected in \cite[Lemma~2.6]{PolstraUSTP}. We record the following consequence.

\begin{lemma}
    \label{lem: symbolic powers contained in splitting ideals}
    Let $R$ be a strongly $F$-regular $F$-finite domain of prime characteristic $p>0$ and Krull dimension $d$. Let $I\subseteq R$ be a nonzero proper ideal and $h$ the maximal height of an associated prime of $R/I$. Then $I^{(hp^e)}\subseteq I_e(I)$ for every $e\geq 1$.
\end{lemma}

\begin{proof}
    By \cite[Lemma~2.6(f)]{PolstraUSTP}, it suffices to check the containment after localizing at the associated primes of $R/I$. The result then follows from \cite[Theorem~13.2.1 and Corollary~13.2.2]{HunekeSwanson2006}.
\end{proof}

Splitting ideals provide a convenient connection between Frobenius splittings, $F$-signature, and $F$-pure thresholds. Let $(R,\fm,k)$ be an $F$-finite local domain of dimension $d$, and let $a_e(R)$ denote the maximal rank of a free $R$-summand of $F^e_*R$. By \cite[Lemma~2.1 and Remark~2.3(1)]{Yao2006}, the normalized splitting numbers satisfy
\[
    \frac{a_e(R)}{\rank_R(F^e_*R)}
    =\frac{\mathcal{L}_R(R/I_e(\fm))}{p^{ed}}.
\]
In particular, $s(R)=\lim_{e\to\infty}\mathcal{L}_R(R/I_e(\fm))/p^{ed}$.

Suppose further that $R$ is $F$-pure, and let $0\not=\fa\subseteq\fm$ be an ideal. For each $e\geq 1$, set $\nu_e(\fa)=\max\{n\geq 0\mid \fa^n\not\subseteq I_e(\fm)\}$, where $\fa^0=R$. The $F$-pure threshold of $\fa$, introduced in \cite{MTW05,TW04}, is $\FPT_R(\fa)=\lim_{e\to\infty}\nu_e(\fa)/p^e$. For $0\not=x\in\fm$, we write $\fpt_R(x)=\FPT_R(xR)$.

If $(R,\fm,k)$ is a local $F$-finite $F$-pure domain, $0\not=x\in \fm$, then $0\leq \nu_e(x)<p^e$ as $x^{p^e}\in \fm^{[p^e]}\subseteq I_e(\fm)$. Moreover, if $\nu_e(x)>0$ for some $e\in\NN$, then $\fpt_R(x)>0$ by \cite[Proposition~2.2(3)]{TakagiWatanabe2004}. Finally, recall that the $\alpha_F$-invariant of an $F$-finite $F$-pure local domain $(R,\fm,k)$ is $\alpha_F(R)=\inf\{\overline{\ord}_{\fm}(x)\fpt_R(x)\mid 0\not=x\in\fm\}$. If $m_e(R)=\max\{n\geq 0\mid I_e(\fm)\subseteq\overline{\fm^n}\}$, then \cite[Proposition~3.7]{TY26} gives
\[
    \alpha_F(R):=\inf\{\overline{\ord}_{\fm}(x)\fpt_R(x)\mid 0\not=x\in\fm\}=\lim_{e\to\infty}\frac{m_e(R)}{p^e}.
\]

\subsection{Models and Reduction to prime characteristic}
Let $(R,\fm,k)$ be a normal local domain essentially of finite type over $\CC$ and set $X=\Spec(R)$. Choose a presentation
\[
  R\cong\bigl(\mathbb C[x_1,\ldots,x_n]/I\bigr)_{\mathfrak p},
\]
where $I\subseteq\mathfrak p$ are ideals of $\mathbb C[x_1,\ldots,x_n]$ and $\mathfrak p$ is prime.  After taking an affine neighborhood of the point, we may assume that the ambient finite-type $\mathbb C$-algebra is itself normal. Taking finitely many coefficients of generators of $I$ and $\mathfrak p$, we
obtain a finitely generated $\mathbb Z$-subalgebra $B\subseteq\mathbb C$, an ideal $I_B\subseteq B[x_1,\ldots,x_n]$, and a prime ideal $\fp_B\supseteq I_B$ such that
\[
  I_B\mathbb C[x_1,\ldots,x_n]=I,
  \qquad
  \mathfrak p_B\mathbb C[x_1,\ldots,x_n]=\mathfrak p.
\]
Set $S_B=B[x_1,\ldots,x_n]/I_B$ and $X_B=\Spec(S_B)$.  Then $X_B\times_{\Spec(B)}\Spec(\CC)$ is identified with $\Spec(\mathbb C[x_1,\ldots,x_n]/I)$, and the prime
$\mathfrak p_BS_B$ spreads out the point corresponding to $\mathfrak p$. We call the data $(B,S_B,\mathfrak p_BS_B)$ a \emph{model} of the local ring $R$.  

By generic freeness and standard spreading-out results, we may replace $B$ by $B[b^{-1}]$ for some $0\ne b\in B$ and assume that the relevant families are flat, that the fibers of $X_B$ are normal and integral, and that the fibers of the closure of the chosen point are integral. For a closed point $\mu\in\Spec(B)$, let $\kappa(\mu)=B/\mu$ and define
\[
  S_\mu=S_B\otimes_B\kappa(\mu),\qquad
  \mathfrak p_\mu=(\mathfrak p_BS_B)S_\mu,\qquad
  R_\mu=(S_\mu)_{\mathfrak p_\mu}.
\]
After shrinking $\Spec(B)$ further, $\mathfrak p_\mu$ is prime for each closed point under consideration.  Since $B$ is a finitely generated $\mathbb Z$-algebra, $\kappa(\mu)$ is a finite field; thus every $R_\mu$ is an $F$-finite local ring of positive characteristic.  

If an effective $\QQ$-divisor $\Delta$ on $X$ is under consideration, we can similarly spread out its finitely many prime components and rational coefficients, obtaining a divisor $\Delta_B$ and reductions $\Delta_\mu$.  After further shrinking the base, the properties needed to discuss the $\mathbb Q$-Cartier divisor $K_{R_\mu}+\Delta_\mu$ are preserved.  Throughout, a statement that holds for \emph{general reductions} means one that holds for every closed point in some dense open subset of $\Spec(B)$.

The local ring $R$ is said to be of \emph{strongly $F$-regular type}
if, for a model as above, there is a dense open subset
$U\subseteq\Spec(B)$ such that $R_\mu$ is strongly $F$-regular for every closed point $\mu\in U$.  This condition is independent of the choice of model.  The corresponding condition for pairs is defined by requiring $(R_\mu,\Delta_\mu)$ to be strongly $F$-regular.  If $K_X+\Delta$ is $\mathbb Q$-Cartier, the fundamental comparison between multiplier and test ideals yields $(X,\Delta)$ is klt if and only if $(R_\mu,\Delta_\mu)$ is strongly $F$-regular for general $\mu$, see \cite{HaraWatanabe2002,Takagi2004}, or \cite[Theorem~2.11]{TY26}.  Following \cite[Definition~4.1]{TY26}, we define the characteristic zero $F$-signature of $R$ by
\[
    s(R)
    := \sup_{\substack{U\subseteq\Spec(B)\\U\text{ dense open}}} \left(\inf_{\substack{\mu\in U\\\mu\text{ closed}}}s(R_\mu)\right).
\]
By \cite[Proposition~4.2]{TY26}, this quantity is independent of the choice of model $(B,S_B,\mathfrak p_BS_B)$.

An application of Main Theorem~\ref{main theorem: F-sig lowerbound via fpt of single element} is that $s(R)>0$ whenever $R$ is a local $\CC$-algebra essentially of finite type and of klt type. In particular, this will establish
\cite[Conjecture~A]{TY26}. To establish this result, we require that the constant $E$ appearing in Main Theorem~\ref{main theorem: F-sig lowerbound via fpt of single element} can be chosen uniformly over the general prime characteristic fibers of a suitable model of $R$. Takagi and Yamaguchi provide a suitable Izumi-Rees type inequality in \cite[Lemma~4.14]{TY26} to bound the constant $E$ appearing in Main Theorem~\ref{main theorem: F-sig lowerbound via fpt of single element} in a reduction to prime characteristic setting.

\begin{proposition}[{\cite[Proof of Lemma~4.14]{TY26}}]
    \label{prop: E in mod p reduction}
    Let $(R,\fm,k)$ be a normal local domain essentially of finite type over $\CC$. Let $R_B$ be a model of $R$ over a finitely generated $\ZZ$-subalgebra $B\subseteq \CC$. There exists an integer $E>0$ so that for a dense open set $U\subseteq \Spec(B)$, every closed point $\mu\in U$, and every $0\not=x\in\fm_\mu$, $e(R_\mu/xR_{\mu})\leq E\overline{\ord}_{\fm_\mu}(x)$.
\end{proposition}
\begin{proof}
    By the final paragraph of the proof of \cite[Lemma~4.14]{TY26}, there exists a constant $C$ and dense open set of closed points $U\subseteq \Spec(B)$ so that for every $\mu\in U$, for each Rees valuation $\nu\in\mathcal{R}_{\fm_{\mu}}$, and every $0\not=x\in\fm_{\mu}$, $\nu(x)\leq C\overline{\ord}_{\fm_{\mu}}(x)$. By \cite[Proposition~3.1]{Lipman1982}, cf \cite[Theorem~3.4 and Corollary~3.5]{Smirnov2020}, we can replace $U$ by a smaller dense open set of closed points and assume for each $\mu\in U$ that $e(R_{\mu}) = e(R_{{\mu'}})$. Let $e_U$ denote the common value of the Hilbert-Samuel multiplicity $e(R_{\mu})$ for $\mu\in U$.

    By the main result of \cite{Rees1961}, for general closed point $\mu$, there are positive integers $d_{\mu,\nu}$ for each $\nu\in\mathcal{R}_{\fm_{\mu}}$ so that
    \[
    e(R_{\mu}/xR_{\mu}) = \sum_{\nu\in \mathcal{R}_{\fm_\mu}}\nu(x)d_{\mu,\nu} \quad \mbox{and} \quad e(R_{\mu}) = \sum_{\nu\in \mathcal{R}_{\fm_\mu}}\nu(\fm_{\mu})d_{\mu,\nu}
    \]
    see \cite[Theorem~2.8 and Corollary 2.9]{Polstra2026} and \cite[Theorem~1.1]{CutkoskyMontano2025} for details. It follows that for every $\mu\in U$, $e(R_{\mu}/xR_{\mu})\leq Ce_U\overline{\ord}_{\fm_{\mu}}(x)$.
\end{proof}

\section{Jacobian Elements, Field Norm Maps, and Lower Bounds of the $F$-signature}
Let $A\subseteq R$ be a finite, generically \'etale extension of excellent noetherian domains, with $A$ regular. We summarize background on relative Jacobian ideals found in \cite{LipmanSathaye1981,Hochster2002}. Choose a presentation $R\cong A[z_1,\ldots,z_c]/P$. Then $P\subseteq A[z_1,\ldots,z_c]$ is a prime ideal of height $c$. If $P=(g_1,\ldots,g_t)$, then the relative Jacobian ideal $\mathcal{J}_{R/A}$ is the ideal of $R$ generated by the images of the $c\times c$ minors of the $t\times c$ matrix $(\frac{\partial g_i}{\partial z_j})$. This ideal is independent of the chosen presentation and is nonzero because $A\subseteq R$ is generically \'etale. To identify a nonzero element of $\mathcal{J}_{R/A}$, observe that $P$ is a prime ideal of height $c$ in the regular ring $A[z_1,\ldots,z_c]$. Thus, we may choose $f_1,\ldots,f_c\in P$ that generate $P$ upon localization at $P$. Since $A\subseteq R$ is generically \'etale, the image $\delta$ of the Jacobian determinant $\det\left(\frac{\partial f_i}{\partial z_j}\right)$ in $R$ is nonzero.

The theory of Jacobian ideals has played an important role in several developments in commutative algebra and algebraic geometry. Of particular relevance to the present article is its connection to the tight closure and test element theories of Hochster and Huneke \cite{HochsterHunekeTC}. If $A$ and $R$ have prime characteristic $p>0$, then every nonzero element of $\mathcal{J}_{R/A}$, including $\delta$, is a completely stable test element of $R$ by \cite[Corollary~(1.5.4)]{HochsterHunekeTC}. Consequently, $R$ is strongly $F$-regular if and only if there exists $e_0\in\NN$ such that $F^{e_0}_*\delta$ generates a free $R$-summand of $F^{e_0}_*R$.

Another important application of Jacobian elements arises in the formal implicit function theorem \cite[Corollary on p.~A.IV.37]{BourbakiAlgebraII}. This implicit function theorem motivates techniques of \cite{Samuel1951Singularites,Samuel1956Algebricite,Hironaka1965Equivalence,CutkoskySrinivasan1993Intrinsic,CutkoskySrinivasan1997Equivalence} and provides criteria for quotients of formal power series rings over a field $k$ to be isomorphic via isomorphisms induced by continuous $k$-algebra automorphisms of the ambient power series ring.

We take inspiration from the argument appearing in \cite[\text{Lemma 6.3}]{CutkoskySrinivasan1993Intrinsic}. First we observe that if $k$ is a field, $x_1,\ldots,x_n$ and $w_1,\ldots,w_n$ are variables, then by binomial expansion, if $a_1,\ldots,a_n\geq 0$
\[
\prod_{i=1}^n(x_i+w_i)^{a_i} = \prod_{i=1}^nx_i^{a_i} + \sum_{i=1}^n a_i\left(x_1^{a_1}\cdots x_{i-1}^{a_{i-1}}x_i^{a_i-1}x_{i+1}^{a_{i+1}}\cdots x_n^{a_n}\right)w_i + H
\]
where $H\in (w_1,\ldots,w_n)^{2}k[x_1,\ldots,x_n,w_1,\ldots,w_n]$.  Consequently, for any $f\in k[[x_1,\ldots,x_n]]$
\[
f(x_1+w_1,\ldots,x_n+w_n) = f + \sum_{i=1}^n \frac{\partial f}{\partial x_i}w_i + H
\]
where $H\in (w_1,\ldots,w_n)^{2}k[[x_1,\ldots,x_n,w_1,\ldots,w_n]]$.

\begin{lemma}[Implicit Function Lemma]
    \label{lem: multivariable calculus lemma}
    Let $k$ be a field, $\underline{x} = x_1,\ldots,x_d$ and $\underline{z} = z_1,\ldots,z_c$ variables,  and $ k[[\underline{x},\underline{z}]]$ the power series ring over $k$. Let $f_1,\ldots,f_c$ be a sequence in $k[[\underline{x},\underline{z}]]$ and $D = \det(\frac{\partial f_i}{\partial z_{j}})_{1\leq i,j\leq c}$. Let $\underline{y} = y_1,\ldots,y_d$ be variables. There exist elements $\lambda_1,\ldots,\lambda_c\in \left(\underline{y}\right)k[[\underline{x},\underline{z},\underline{y}]]$ and a continuous $k$-linear algebra homomorphism 
    \[
    \Psi:k[[\underline{x},\underline{z}]]\to k[[\underline{x},\underline{z}, \underline{y}]]
    \]
    such that
    \begin{itemize}
        \item $\Psi(x_i) = x_i+D^2y_i$ for all $1\leq i\leq d$;
        \item $\Psi(z_i) = z_i + \lambda_i$ for all $1\leq i\leq c$;
        \item $\Psi(f_i) = f_i$ for all $1\leq i\leq c$.
    \end{itemize}
\end{lemma}

\begin{proof}
    Let $\underline{u} = u_1,\ldots,u_d$ and $\underline{w} = w_1,\ldots,w_c$ denote variables. Consider the expansion of $f_i(\underline{x}+\underline{u},\underline{z}+\underline{w})$. For each $1\leq i\leq c$ we can write
    \[
    f_i(\underline{x}+\underline{u},\underline{z}+\underline{w}) = f_i(\underline{x},\underline{z}) + \sum_{j = 1}^d\frac{\partial f_i}{\partial x_j}u_j + \sum_{j=1}^c\frac{\partial f_i}{\partial z_j}w_j + H_i(\underline{u},\underline{w})
    \]
    where $H_i(\underline{u},\underline{w})\in (\underline{u},\underline{w})^2k[[\underline{x},\underline{z},\underline{u},\underline{w}]]$.

    Let $\underline{y} = y_1,\ldots,y_d$ and $\underline{v} = v_1,\ldots,v_c$ denote variables and let $\underline{\omega} = \omega_1,\ldots,\omega_c$ with $\omega_i$ the $i$th entry of $\adj(\frac{\partial f_i}{\partial z_j})[v_1,\ldots,v_c]^T$. Consider the substitution $\underline{u} = D^2\underline{y}$ and $\underline{w} = D \underline{\omega}$. Then $\sum_{j=1}^c\frac{\partial f_i}{\partial z_j}\omega_j = D v_i$ and
    \begin{align*}
        f_i(\underline{x}+D^2\underline{y},\underline{z}+D \underline{\omega}) &= f_i(\underline{x},\underline{z}) + \sum_{j = 1}^d\frac{\partial f_i}{\partial x_j}D^2y_j + D^2 v_i + H_i(D^2\underline{y},D \underline{\omega})
    \end{align*}
    with $H_i(D^2\underline{y},D \underline{\omega})\in (D^2\underline{y},D\underline{v})^2k[[\underline{x},\underline{z},\underline{y},\underline{v}]]$. The sum $\sum_{j = 1}^d\frac{\partial f_i}{\partial x_j}D^2y_j + D^2 v_i + H_i(D^2\underline{y},D \underline{\omega})$ is divisible by $D^2$ and is of the form
    \[
    \sum_{j = 1}^d\frac{\partial f_i}{\partial x_j}D^2y_j + D^2 v_i + H_i(D^2\underline{y},D \underline{\omega}) = D^2\left(v_i+\sum_{j=1}^d\frac{\partial f_i}{\partial x_j}y_j+Q_i(\underline{y},\underline{v})\right)
    \]
    with $Q_i(\underline{y},\underline{v})\in (\underline{y},\underline{v})^2
    k[[\underline{x},\underline{z},\underline{y},\underline{v}]]$. For each $1\leq i\leq c$ let $g_i(\underline{y},\underline{v}) = v_i+\sum_{j=1}^d\frac{\partial f_i}{\partial x_j}y_j+Q_i(\underline{y},\underline{v})$ so that for each $1\leq i\leq c$,
    \begin{equation}
        \label{eq: subsitution for f}
        f_i(\underline{x}+D^2\underline{y},\underline{z}+D \underline{\omega}) = f_i(\underline{x},\underline{z}) + D^2g_i(\underline{y},\underline{v}).
    \end{equation}
    
    The matrix of partial derivatives $(\frac{\partial g_i}{\partial v_j})_{1\leq i,j\leq c}$ evaluated at $\underline{y} = \underline{v} = 0$ is the identity matrix. Hence the constant term of $\det(\frac{\partial g_i}{\partial v_j})_{1\leq i,j\leq c}$ is invertible in $k[[\underline{x},\underline{z}]]$ and by the formal implicit function theorem, \cite[Corollary on p.~A.IV.37]{BourbakiAlgebraII}, there is a unique tuple of formal power series $\underline{\gamma} = \gamma_1,\ldots,\gamma_c\in (\underline{y})k[[\underline{x},\underline{z},\underline{y}]]$ so that $g_i(\underline{y},\underline{\gamma}) = 0$ for all $1\leq i \leq c$.

   Let $\underline{\lambda}$ denote the sequence in $(\underline{y})k[[\underline{x},\underline{z},\underline{y}]]$ obtained from the sequence $D\underline{\omega}$ upon replacing $v_i$ by $\gamma_i$. Then the substitution $v_i = \gamma_i$ has the property for each $1\leq i\leq c$
   \begin{align*}
   f_i(\underline{x}+D^2\underline{y},\underline{z}+\underline{\lambda}) &= f_i(\underline{x},\underline{z}) + D^2g_i(\underline{y},\underline{\gamma})\\
   &=f_i(\underline{x},\underline{z}).
   \end{align*}
   Therefore the continuous $k$-algebra homomorphism $\Psi: k[[\underline{x},\underline{z}]] \to k[[\underline{x},\underline{z},\underline{y}]]$ defined by $\Psi(x_i) = x_i+D^2y_i$ for each $1\leq i\leq d$ and $\Psi(z_i) = z_i+\lambda_i$ for each $1\leq i\leq c$ also satisfies $\Psi(f_i)= f_i$ for all $1\leq i\leq c$.
\end{proof}

\begin{corollary}
    \label{cor: mapping to discriminant power}
    Let $k$ be a field of prime characteristic $p>0$, $\underline{x} = x_1,\ldots,x_d$ and $\underline{z} = z_1,\ldots,z_c$ variables, $S = k[[\underline{x},\underline{z}]]$, $P\in\Spec(S)$ a prime ideal of height $c$ so that $A=k[[\underline{x}]]\subseteq R:=S/P$ is finite and generically \'{e}tale. Denote by $\fm_A,\fm_R,\fm_S$ the maximal ideals of $A, R,$ and $S$ respectively. Let $f_1,\ldots,f_c\in P$ be local generators of $P$ and $D=\det(\frac{\partial f_i}{\partial z_j})_{1\leq i,j\leq c}$. Let $\delta$ denote the image of $D$ in the quotient ring $R$. If $0\not=g\in A$ and $p^e>\ord_{\fm_A}(g)$, then there exists $\Phi\in\Hom_R(F^e_*R,F^e_*R)$ so that $\Phi(F^e_*g) = F^e_*\delta^{2\ord_{\fm_A}(g)}$.
\end{corollary}

\begin{proof}
    The existence of $\Phi\in \Hom_R(F^e_*R,F^e_*R)$ so that $\Phi(F^e_*g) = F^e_*\delta^{2\ord_{\fm_A}(g)}$ is equivalent to finding $\Phi'\in \Hom_{R^{p^e}}(R,R)$ so that $\Phi'(g) = \delta^{2\ord_{\fm_A}(g)}$. By \cref{lem: multivariable calculus lemma}, there exists $\lambda_1,\ldots, \lambda_c\in (\underline{y})k[[\underline{x},\underline{z},\underline{y}]]$ and a $k$-linear algebra map $\Psi:k[[\underline{x},\underline{z}]]\to k[[\underline{x},\underline{z},\underline{y}]]$ so that 
    \begin{itemize}
        \item $\Psi(x_i) = x_i+D^2y_i$ for all $1\leq i\leq d$;
        \item $\Psi(z_i) = z_i + \lambda_i$ for all $1\leq i\leq c$;
        \item $\Psi(f_i) = f_i$ for all $1\leq i\leq c$.
    \end{itemize}
    Observe that $\Psi(P)\subseteq Pk[[\underline{x},\underline{z},\underline{y}]]$ as the containment can be checked upon localization at $Pk[[\underline{x},\underline{z},\underline{y}]]$ where $f_1,\ldots,f_c$ remain local generators. Therefore $\Psi$ induces $k$-algebra map $\Psi_R: R\to R[[\underline{y}]]$ so that 
    \begin{itemize}
        \item $\Psi_R(x_i) = x_i+\delta^2y_i$ for all $1\leq i\leq d$;
        \item $\Psi_R(z_i) = z_i + \lambda_i$ for all $1\leq i\leq c$.
    \end{itemize}
    Consider the induced map $\overline{\Psi}_R: R \to R[[\underline{y}]]/(\underline{y})^{\ord_{\fm_A}(g)+1}$. We note that $R[[\underline{y}]]/(\underline{y})^{\ord_{\fm_A}(g)+1}$ is a free $R$-module with a basis of monomials in the $y$'s of degree at most $\ord_{\fm_A}(g)$. 
    
    Since $y_i^{p^e},\lambda_j^{p^e}\in (\underline{y})^{\ord_{\fm_A}(g)+1}$ for each $1\leq i\leq d$ and $1\leq j\leq c$, $\overline{\Psi}_R(x_i^{p^e}) = (x_i + \delta^2 y_i)^{p^e} = x_i^{p^e}$ and $\overline{\Psi}_R(z_j^{p^e}) = (z_j+\lambda_j)^{p^e} = z_j^{p^e}$. Therefore $\overline{\Psi}_R: R\to R[\underline{y}]/(\underline{y})^{\ord_{\fm_A}(g)+1}$ is $R^{p^e}$-linear.

    We choose $\beta \in \bN^d$ with $|\beta|= \ord_{\fm_A}(g)$ such that the coefficient of the monomial $x^{\beta}$ in the homogeneous component of degree $\ord_{\fm_A}(g)$ of $g$ is nonzero. Let $\pi$ be the projection map from $R[[\underline{y}]]/(\underline{y})^{\ord_{\fm_A}(g)+1}$ onto the monomial $y^{\beta}$. It is easy to check that $(\pi \circ \overline{\Psi}_R)(g)$ is a unit multiple of $\delta^{2\ord_{\fm_A}(g)}$ in $R$. Therefore there exists $\Phi'\in\Hom_{R^{p^e}}(R,R)$ so that $\Phi'(g) = \delta^{2\ord_{\fm_A}(g)}$.
    \end{proof}
    

Let $A\subseteq R$ be a finite inclusion of excellent domains with respective fraction fields $K$ and $L$. We denote by $N_{L/K}:L\to K$ the field norm map. For $x\in L$, let $T_x:L\to L$ be the multiplication by $x$ map viewed as an element of $\Hom_K(L,L)$, then $N_{L/K}(x)$ is the determinant of $T_x$. If $x\in R$, then $N_{L/K}(x)\in K$ is integral over $A$. Hence if $A$ is normal then $N_{L/K}(R)\subseteq A$ and we then denote by $N_{R/A}:R\to A$ the field norm map of $K\subseteq L$ restricted to $R$.

\begin{lemma}
    \label{lem: Field norm map lemma}
    Let $(A,\fm_A,\kappa(A))\subseteq (R,\fm_R,\kappa(R))$ be a finite inclusion of excellent local domains so that $A$ is normal and $R$ is analytically irreducible. Let $N_{R/A}:R\to A$ be the field norm map restricted to $R$.
    \begin{enumerate}
        \item For every $x\in R$, $N_{R/A}(x)\in xR\cap A$.
        \item If $0\not=x\in \fm_R$, then 
        \[
        e\left(\frac{A}{N_{R/A}(x)A}\right) = [\kappa(R):\kappa(A)] \cdot e_{\frac{(\fm_A,x)R}{xR}}\left(\frac{R}{x R}\right).
        \]
        In particular, if $A$ is a regular local ring then
        \[
        \ord_{\fm_A}(N_{R/A}(x)) = [\kappa(R):\kappa(A)] \cdot e_{\frac{(\fm_A,x)R}{xR}}\left(\frac{R}{x R}\right).
        \]
    \end{enumerate}
\end{lemma}

\begin{proof}
    Let $x\in R$, $T_x\in \Hom_A(R,R)$ the multiplication by $x$ map. Let $K,L$ denote the fraction fields of $A$ and $R$ respectively. Let $p(t)\in K[t]$ be the characteristic polynomial of $T_x$. The coefficients of $p(t)$ are elements of $A$ and the constant coefficient of $p(t)$ is $\pm N_{R/A}(x)$. Hence there exists a polynomial $q(t)\in A[t]$ so that $T_x\circ q(T_x) = N_{R/A}(x)$. Evaluating the equal functions at $1\in R$, 
    \[
    N_{R/A}(x) = (T_x\circ q(T_x))(1) = xq(T_x)(1).
    \]
    Therefore $N_{R/A}(x)\in xR\cap A$.

    Since for every finite-length $R$-module $M$, we have $\mathcal{L}_A(M)= [\kappa(R):\kappa(A)]\mathcal{L}_R(M)$, we know that 
    \[
        e_{\fm_A}\left(\frac{R}{xR}\right) = [\kappa(R):\kappa(A)] e_{\fm_A R}\left(\frac{R}{x R}\right) = [\kappa(R):\kappa(A)]  e_{\frac{(\fm_A,x)R}{xR}}\left(\frac{R}{x R}\right).
    \]
    Let $\mathfrak p\subseteq A$ be a height-one prime ideal. Since $A$ is normal, $A_{\mathfrak{p}}$ is a DVR. Furthermore, since $R$ is finite and torsion-free over $A$, $R_{\mathfrak{p}} = R \otimes_A A_{\mathfrak{p}}$ is a finite free $A_{\mathfrak{p}}$-module. Let $\pi$ be the uniformizer of $A_{\mathfrak{p}}$. By Smith normal form, we can represent $T_x$, after changing bases, by a diagonal matrix ${\rm diag}(\pi^{n_1},\cdots,\pi^{n_r})$, up to multiplication by units. Thus,
    \[
     R_{\mathfrak{p}}/xR_{\mathfrak{p}} \cong \oplus_{i=1}^{r}A_{\mathfrak{p}}/\pi^{n_i}A_{\mathfrak{p}}.
    \]
    Taking lengths gives
    \[
    \mathcal{L}_{A_{\mathfrak{p}}}(R_{\mathfrak{p}}/xR_{\mathfrak{p}}) = \sum_{i=1}^{r}n_i = v_{\mathfrak{p}}({\rm det}T_x) = v_{\mathfrak{p}}(N_{R/A}(x)) = \mathcal{L}_{A_{\mathfrak{p}}}((A/N_{R/A}(x)A)_{\mathfrak{p}}).
    \]
    The equality of multiplicities is then an application of the associativity formula for multiplicity.
\end{proof}

\begin{theorem}
    \label{thm: alpha F-invariant lower bound by fpt of Jacobian element}
    Let $(R,\fm,k)$ be a complete $F$-finite local strongly $F$-regular domain of prime characteristic $p>0$ and Krull dimension $d$ that is not regular. Let $A = k[[x_1,\ldots,x_d]]\subseteq R$ be a regular local subring of $R$ so that $A\subseteq R$ is finite and generically \'etale. Fix a presentation $R\cong A[z_1,\ldots,z_c]/P$. Let $f_1,\ldots,f_c\in P$ generate $PA[z_1,\ldots,z_c]_P$ and $D=\det(\frac{\partial f_i}{\partial z_j})_{1\leq i,j\leq c}$. Let $\delta$ denote the image of $D$ in the quotient ring $R$. Let $E$ be a constant so that for all $0\not=x\in \fm_R$, 
        \[
        e_{\frac{(\fm_A,x)R}{xR}}\left(\frac{R}{xR}\right)\leq E\overline{\ord}_{\fm_R}(x).
        \]
        Then
        \[
        \alpha_F(R)\geq \frac{\fpt_R(\delta)}{2E}\quad \mbox{and}\quad s(R) \geq \frac{e(R)}{d!}\left(\frac{\fpt_R(\delta)}{2E}\right)^d.
        \]
\end{theorem}

\begin{proof}
    Denote by $\fm_A,\fm_R$ the maximal ideals of $A$ and $R$ respectively. It suffices to prove the claimed lower bound of $\alpha_F(R)$, see \cite[Proposition~3.9 (1)]{TY26}. Because we are assuming $R$ is not regular, $\delta\in \fm_R$. The claims of the theorem are trivial if $\fpt_R(\delta) = 0$. We therefore may assume $\fpt_R(\delta)>0$, implying that $R$ is strongly $F$-regular. Let $m_e(R) = \max\{n\in\NN\mid I_e(\fm)\subseteq \overline{\fm^n}\}$ so that $\alpha_F(R) = \lim_{e\to\infty}\frac{m_e(R)}{p^e}$.
    
    If $p^e \leq \lceil E (m_e(R)+1)\rceil$ for infinitely many $e\in\NN$, then the desired inequality is trivial since
    \[
    \frac{\fpt_R(\delta)}{2E}\leq \frac{1}{E} \leq \lim_{e\to \infty}\frac{\lceil E (m_e(R)+1)\rceil}{Ep^e} = \alpha_F(R)
    \]
    where the limit is taken along the infinitely many choices of $e$ so that $p^e \leq \lceil E (m_e(R)+1)\rceil$. We therefore may assume for $e\gg 0$ that $p^e>\lceil E (m_e(R)+1)\rceil$. Let $e\gg 0$ and choose $g\in I_e(\fm)$ so that $m_e(R)\leq \overline{\ord}_{\fm_R}(g) < m_e(R) + 1$. By \cref{lem: Field norm map lemma},  
    \[
    N_{R/A}(g)\in gR\cap A\subseteq I_e(\fm)\quad \mbox{and}\quad\ord_{\fm_A}(N_{R/A}(g)) = e_{\frac{(\fm_A,x)R}{(x)}}\left(\frac{R}{gR}\right).
    \]
    By assumption,
    \[
    \ord_{\fm_A}(N_{R/A}(g)) = e_{\frac{(\fm_A,g)R}{gR}}\left(\frac{R}{g R}\right)\leq E\overline{\ord}_{\fm_R}(g) < E(m_e(R)+1)\leq \lceil E (m_e(R)+1)\rceil<p^e.
    \]
    By \cref{cor: mapping to discriminant power}, there exists $\Psi\in \Hom_R(F^e_*R,F^e_*R)$ so that $\Psi(F^e_*N_{R/A}(g)) = F^e_*\delta^{2\lceil E (m_e(R)+1)\rceil}$.  Therefore $\delta^{2\lceil E(m_e(R)+1)\rceil}\in I_e(\fm)$, else if $\varphi\in\Hom_R(F^e_*R,R)$ is so that $\varphi(F^e_*\delta^{2\lceil E (m_e(R)+1)\rceil}) = 1$, then $\varphi\circ\Psi\in \Hom_R(F^e_*R,R)$ is such that $(\varphi\circ\Psi)(F^e_*N_{R/A}(g)) = 1$, contradicting $N_{R/A}(g)\in I_e(\fm)$.

    Let $\nu_e(\delta) = \max\{t\in\NN\mid \delta^t\not\in I_e(\fm)\}$. Then 
    \[
    \nu_e(\delta)<2\lceil E (m_e(R)+1)\rceil \quad\mbox{and}\quad \fpt(\delta) = \lim_{e\to\infty}\frac{\nu_e(\delta)}{p^e}\leq \lim_{e\to \infty}\frac{2\lceil E (m_e(R)+1)\rceil}{p^e} = 2E\alpha_F(R).
    \]
    Dividing by $2E$ gives the claimed lower bound of $\alpha_F(R)$. 
\end{proof}

\begin{remark}
    With the notation of \cref{thm: alpha F-invariant lower bound by fpt of Jacobian element}, if $\overline{\fm_AR} = \fm_R$ and $\Gr_{\fm}(R)$ is a normal domain, e.g., $A\subseteq R$ is obtained by completing a standard graded normal domain over $k$ and $x_1,\ldots,x_d$ are obtained from degree $1$ elements, then the constant $E$ appearing in \cref{thm: alpha F-invariant lower bound by fpt of Jacobian element} can be taken to be $e(R)$. Indeed, for every $0\not=x\in\fm$, $\ord_{\fm}(x) = \overline{\ord}_{\fm}(x)$ and $e_{\frac{(\fm_A,x)R}{xR}}(R/xR) =e(R/xR)= e(R)\ord_{\fm}(x)$.
\end{remark}

\begin{corollary}
    \label{cor: positive lim F-signature}
    Let $(R,\fm,k)$ be a local normal domain essentially of finite type over $\CC$. If $R$ has klt type then $s(R) >0.$
\end{corollary}

\begin{proof}
    Let $T$ be a finitely generated $\CC$-algebra so that $R\cong T_{\fp}$ for some $\fp\in \Spec(T)$. We may assume the finitely generated $\CC$-algebra $T$ is a normal domain. Let $\fq$ be a closed point of $T$ in a klt neighborhood of $\fp$ of the pair $(T,\Delta)$. since $F$-signature does not decrease upon localizations, see \cite{Polstra2018,PolstraTucker2018}, it suffices to show the $F$-signature of a closed point in the klt locus of a boundary pair of a finitely generated normal domain has positive $F$-signature. We therefore assume $\fp$ to be a closed point of $\Spec(T)$ and set $d = \dim(R)$.
    
     Let $A = \mathbb{C}[x_1,\ldots,x_d]$ be a choice of polynomial subring of $T$ so that $x_1,\ldots,x_d\in \fp$ and $A\subseteq T$ is finite. Let $z_1,\ldots,z_c\in \fp$ generate $T$ as an $A$-algebra. Abusing notation, let $\CC[x_1,\ldots,x_d,z_1,\ldots,z_c]$ be a polynomial ring mapping onto $T$ so that we may identify $T\cong \CC[x_1,\ldots,x_d,z_1,\ldots,z_c]/P$ for some prime $P$ of height $c$. Let $f_1,\ldots,f_c\in P$ generically generate $P$, $D = \det(\frac{\partial f_i}{\partial z_j})_{1\leq i,j \leq c}$, and $\delta$ the (nonzero) image of $D$ in the quotient ring $T$ since $T$ generically \'{e}tale over $\CC[x_1,\ldots,x_d]$ by characteristic considerations. We may assume $\delta\in \fp$, else $R$ is regular and a general fiber of a model of $R$ is regular.

    We choose a model of $R$ over a finitely generated $\ZZ$-algebra $B\subseteq \CC$. Completion of the finite normalization decomposes into completed local factors, one of which is $\widehat{R_\mu}$. Thus for a dense open set $U\subseteq \Spec(B)$, generic \'{e}taleness and the chosen generic generators are preserved at the closed points of $U$. Thus for every sufficiently general closed point $\mu\in \Spec(B)$, the completed local ring $\widehat{R_\mu}$ has an induced presentation
    \[
    \widehat{R_\mu} \cong \kappa(\mu)[[\underline{x},\underline{z}]]/\widehat P_\mu
    \]
    and is finite and generically \'{e}tale over
    \[
    \widehat A_\mu= \kappa(\mu)[[\underline{x}]].
    \]
    Here $\widehat P_\mu$ is prime, the images of the $f_i$ generate it at its generic point, and the corresponding Jacobian determinant is $\delta_\mu$.

    By \cref{prop: E in mod p reduction}, there exists a constant $E_0$ and a model of $R$ over a finitely generated $\ZZ$-algebra $B\subseteq \CC$, so that for any general closed point $\mu \in \Spec(B)$, $e(R_\mu/xR_{\mu})\leq E_0\overline{\ord}_{\fm_\mu}(x)$. Choose an integer $a\geq 1$ such that $\overline{\fm^a}\subseteq (x_1,\ldots,x_d)R$. After shrinking $U$, this containment specializes to $\overline{\fm_\mu^a}\subseteq (x_1,\ldots,x_d)R_\mu$. Consequently, for each general closed point $\mu$ and $0\not=x\in \fm_{\mu}$, 
    \[
    e_{(x_1,\ldots,x_d)R_\mu}(R_\mu/xR_\mu) \leq a^{d-1}e_{\fm_\mu}(R_\mu/xR_\mu) \leq a^{d-1}E_0\overline{\ord}_{\fm_\mu}(x).
    \]
    Set $E=a^{d-1}E_0$. Both $F$-signature and $F$-pure threshold are constant under completion. Therefore by \cref{thm: alpha F-invariant lower bound by fpt of Jacobian element}, if $\mu \in \Spec(B)$ is a general closed point,
    \[
    s(R_{\mu}) \geq \frac{e(R_{\mu})}{d!}\left(\frac{\fpt_{R_{\mu}}(\delta)}{2E}\right)^d \geq \frac{1}{d!}\left(\frac{\fpt_{R_{\mu}}(\delta)}{2E}\right)^d
    \]

    Choose an effective $\QQ$-Weil divisor $\Delta$ on $X=\Spec(R)$ such that $K_X+\Delta$ is $\QQ$-Cartier and $(X,\Delta)$ is klt. Let $\varepsilon>0$ be a rational number such that $(X,\Delta+\varepsilon\operatorname{div}_X(\delta))$ is klt. After enlarging $B$ to spread out data, \cite[Corollary~3.4]{Takagi2004}, cf \cite[Corollary~6.6]{CEMS18}, implies that there exists a dense open subset $U\subseteq\Spec(B)$ such that $\left(R_\mu,\Delta_\mu+ \varepsilon\operatorname{div}(\delta_\mu)\right)$ is strongly $F$-regular for every closed point $\mu\in U$. Consequently, $\fpt_{R_\mu}(\delta)\geq \fpt_{(R_\mu,\Delta_\mu)}(\delta)\geq \varepsilon$ for every closed point $\mu\in U$. Therefore, 
    \[
    s(R) = \sup_{\substack{V\subseteq\Spec(B)\\ V\text{ dense open}}}\left(\inf_{\substack{\mu\in V\\ \mu\text{ closed}}}s(R_\mu)\right)\geq\inf_{\substack{\mu\in U\\ \mu\text{ closed}}}s(R_\mu)\geq\frac{1}{d!}\left(\frac{\varepsilon}{2E}\right)^d.
    \]
    Allowing $\varepsilon$ to approach $\sup\{\varepsilon\in \QQ_{>0}\mid (R,\Delta + \varepsilon \Div(\delta))\mbox{ is klt}\} = \lct_{(R,\Delta)}(\delta)$, we obtain
    \[
    s(R)\geq \frac{1}{d!}\left(\frac{\lct_{(R,\Delta)}(\delta)}{2E}\right)^d>0 \qedhere
    \]
\end{proof}

Let $X$ be normal irreducible quasi projective variety over $\CC$, $x\in X$ a closed point, and $\Delta\geq 0$ a boundary $\QQ$-divisor. If $x$ is in the klt locus of $(X,\Delta)$ and if $\mathcal{O}_{X,x}$ agrees with the localization of a standard graded normal domain over $\CC$ at its homogeneous ideal, then the following theorem bounds from below the $F$-signature of $\mathcal{O}_{X,x}$ by a constant depending only upon the dimension and embedding dimension of $\mathcal{O}_{X,x}$. We refer the reader to \cite[Section~2.2]{LiXu2020} for background on normalized volumes of klt singularities, including the definitions of log discrepancy $A_{(X,\Delta)}(\nu)$, the volume of a valuation $\vol_{X,x}(\nu)$, the normalized volume of a valuation $\widehat{\vol}_{(X,\Delta)}(\nu)$, and the normalized volume $\widehat{\vol}(x,X,\Delta)$.

\begin{theorem}
    \label{thm: lower bound of F-sig via normalized volume and embedding dimension}
    Let $S = \CC[x_1,\ldots,x_{d+c}]$ be a standard graded polynomial ring, $\fm$ the homogeneous maximal ideal, $0\not=P\subseteq \fm^2$ a homogeneous prime ideal of height $c$, $R = (S/P)_{\fm}$, $X = \Spec(R)$, $x\in X$ the closed point, and $\Delta\geq 0$ a boundary $\QQ$-divisor so that $(X,\Delta)$ is klt. There exists positive constant $C(d,c)$, depending only upon the dimension and embedding dimension of $R$, such that
    \[
    s(R)\geq C(d,c)\left(\widehat{\vol}(x,X,\Delta)\right)^d.
    \]
\end{theorem}

\begin{proof}
    The localization of the $\NN$-graded ring $S/P$ at its homogeneous maximal ideal has at worst rational singularities since $(X,\Delta)$ is klt. Therefore $a(S/P)<0$ by \cite[Theorem~5.1]{KuranoSatoSinghWatanabe2009} where $a(S/P)$ denote the $a$-invariant of $S/P$. The standard graded Cohen-Macaulay domain $S/P$ has regularity $0\leq d+a(S/P)\leq d-1$. Therefore $P$ can be generated by elements of degree at most $d$, see \cite[Section~4a]{Eisenbud2005}. Choose a Noether normalization $A\subseteq S/P$ formed on $d$ general degree $1$ elements and let $z_1,\ldots,z_c$ be degree $1$ elements generating $S/P$ as an $A$-algebra. We may choose $c$ homogeneous minimal generators $f_1,\ldots,f_c$ of $P$ of degree at most $d$ that generate $P$ upon localization at $P$. The resulting Jacobian element $\delta\in S/P$ has degree at most $c(d-1)$.

    Let $E$ be the exceptional divisor of the blowup of $X$ at $x$. The valuation $\nu_E$ of $E$ is the degree valuation of $S/P$. Moreover, for each $0\not=f\in \fm$, $e(R/fR) = \nu_E(f)e(R) = \ord_{\fm}(f)e(R)$. By the proof of \cref{cor: positive lim F-signature}, 
    \[
    s(R) \geq \frac{1}{d!}\left(\frac{\lct_{(X,\Delta)}(\delta)}{2e(R)}\right)^d.
    \]
     By \cite[Lemma~3.4]{Zhuang24}, there exists constant $C_0>0$, depending only on dimension, such that 
    \begin{align*}
    \lct_{(R,\Delta)}(\delta)&\geq C_0 \frac{A_{(X,\Delta)}(\nu_E)}{\widehat{\vol}_{(X,\Delta)}(\nu_E)\nu_E(\delta)}\widehat{\vol}(x,X,\Delta)\\
    &=\frac{C_0}{(A_{(X,\Delta)}(\nu_E))^{d-1}\vol_{X,x}(\nu_E) \nu_E(\delta)}\widehat{\vol}(x,X,\Delta).
    \end{align*}
    The lower bounds of $F$-signature and of $\lct_{(R,\Delta)}(\delta)$ provided above imply that it is enough to find upper bounds of $e(R)$, $A_{(X,\Delta)}(\nu_E)$, $\vol_{X,x}(\nu)$, and $\nu_E(\delta)$ depending only on $d$ and $c$.

    Because $R$ has at worst rational singularities, $e(R)\leq {{d+c-1}\choose{d-1}}$, see \cite[Theorem~3.1]{HunekeWatanabe}. Since $\nu_E$ is the degree valuation, 
    \[
    \vol_{X,x}(\nu_E) = e(R)\leq {{d+c-1}\choose{d-1}}\quad\mbox{and}\quad \nu_E(\delta) = \deg_{\fm}(\delta)\leq c(d-1).
    \]
    To determine discrepancy value $A_{(X,\Delta)}(\nu_E)$, we have $-d<a(S/P)<0$ and as a graded module over $S$, the graded canonical module $\omega_{S/P}$ has a homogeneous element $\eta$ of degree $-a(S/P)$ with $0<-a(S/P)<d$. Let $K_X$ be the effective canonical divisor of $X$ corresponding to $\eta$ and let $\pi:X'\to X$ be the blowup of $\fm$. Then the pullback $\pi^*(K_X+\Delta)$ is effective and
    \begin{align*}
    A_{(X,\Delta)}(\nu_E)&=1+\ord_E(K_{X'}-\pi^*(K_X+\Delta)) \\
    &=1+(-a(S/P)-1) - \ord_E(\pi^*(K_X+\Delta))\\
    &\leq -a(S/P)\\
    &\leq d.\qedhere
    \end{align*}
\end{proof}

\section{Multiplicity Bounds and Ideal Containments}\label{section: effective bounds}

The purpose of this section is to establish effective bounds on multiplicities and ideal containments in terms of $F$-signature. The main technical ingredient is a valuation-theoretic estimate for Frobenius splitting ideals. Specifically, we show that if $x\in I_e(I)$ and $\nu$ is a discrete valuation of a strongly $F$-regular local domain $R$, then $\nu(x)\geq p^es(R)\nu(I)$. The proof uses the free summands of $F^e_*R$ and a determinant computation after extending scalars to the valuation ring.

Combining this estimate with the containment of symbolic powers in splitting ideals and the Rees valuative criterion for integral closure, we obtain explicit uniform symbolic power containments for strongly $F$-regular rings and, through reduction to prime characteristic, for klt singularities. We also establish the multiplicity bound $e(R)\leq d!/s(R)^{d-1}$, whose exponent is optimal by the Veronese examples, and consequently obtain a bound on embedding dimension. Further applications include inequalities relating $F$-signature to the $\alpha_F$-invariant, $F$-pure thresholds, and valuations. In the $\QQ$-Gorenstein setting, these yield an effective Izumi-type inequality involving log discrepancies.

\begin{theorem}
    \label{thm: splitting ideals and valuation rings}
    Let $(R,\fm,k)$ be an $F$-finite strongly $F$-regular local ring of prime characteristic $p>0$ and $K$ the fraction field of $R$. Let $R\subseteq V\subseteq K$ be a discrete valuation ring of $R$ with discrete valuation $\nu$. Then for every ideal $I\subseteq R$, every $e\geq 1$, $e_1\geq 0$, and $x\in I_e(I)$,
    \[
    \nu(x)\geq \frac{p^ea_{e+e_1}(R)}{\rank_R(F^{e+e_1}_*R)}\nu(I).
    \]
    In particular, letting $e_1\to \infty$,
    \[
    \nu(x)\geq p^es(R)\nu(I).
    \]
\end{theorem}

\begin{proof}
    We decompose $F^e_*R\cong R^{\oplus a_e(R)}\oplus N_e$ so that $N_e$ has no free $R$-summands. Let $T_{F^e_*x}:F^e_*R\to F^e_*R$ denote the multiplication by $F^e_*x$ map, i.e., $T_{F^e_*f}(F^e_*x) = F^e_*fx$. Identify $T_{F^e_*x}\in \Hom_R(F^e_*R,F^e_*R)$ as a matrix of $R$-linear maps
    \[
    T_{F^e_*x}: R^{\oplus a_e(R)} \oplus N_e \xrightarrow{\begin{pmatrix}A_{1,1} & A_{1,2} \\ {A_{2,1}} & A_{2,2}\end{pmatrix}} R^{\oplus a_e(R)} \oplus N_e
    \]
    so that $A_{1,1}\in \Hom_R(R^{a_e(R)}, R^{a_e(R)}), A_{2,1}\in \Hom_R(R^{a_e(R)}, N_e), A_{1,2}\in \Hom_R(N_e, R^{a_e(R)})$ and $A_{2,2}\in \Hom_R(N_e, N_e)$. 

    We are assuming $x\in I_e(I)$. Therefore the images of $A_{1,1}$ and $A_{1,2}$ are contained in $I^{\oplus a_e(R)}\subseteq R^{\oplus a_e(R)}$. Next observe $\left(T_{F^e_*x}\right)^{p^e} = T_{F^e_*x^{p^e}} = xT_{F^e_*1} = x\id_{F^e_*R}$, implying
    \[
    \begin{pmatrix}A_{1,1} & A_{1,2} \\ {A_{2,1}} & A_{2,2}\end{pmatrix}^{p^e} = x\begin{pmatrix}\id_{R^{a_e(R)}} & 0\\ 0 & \id_{N_e}\end{pmatrix}.
    \]

    Let $T_{F^e_*x,V}$ denote the induced $V$-linear map obtained from $T_{F^e_*f}$ by base change to $V$,
    \[
    T_{F^e_*x,V}: F^e_*R\otimes_R V\to F^e_*R\otimes_R V.
    \]
    The image of a torsion element of $F_*R\otimes_R V$ under multiplication by $F^e_*x$ remains a torsion element of $F^e_*R\otimes_R V$. Therefore if $\mathcal{T}(F^e_*R\otimes_R V)$ is the torsion submodule of $F^e_*R\otimes_RV$ then $T_{F^e_*x,V}$ induces a well-defined map on the free $V$-module 
    \[
    \frac{F^e_*R\otimes_R V}{\mathcal{T}(F^e_*R\otimes_R V)}\cong V^{\oplus a_e(R)}\oplus \frac{N_e\otimes_RV}{\mathcal{T}(N_e\otimes_RV)}\cong V^{\oplus a_e(R)}\oplus V^{\rank_R(F^e_*R)-a_e(R)}.
    \]
    Let $A_{1,1,V},A_{2,1,V},A_{1,2,V},$ and $A_{2,2,V}$ denote the respective induced $V$-linear homomorphisms obtained by base change of $A_{1,1},A_{2,1}, A_{1,2},$ and $A_{2,2}$ after modding out by the torsion submodules of the domain and codomains of the respective maps.

    The images of $A_{1,1}$ and $A_{1,2}$ contained in $I^{\oplus a_e(R)}\subseteq R^{\oplus a_e(R)}$ imply the images of $A_{1,1,V}$ and $A_{1,2,V}$ are contained in $IV^{\oplus a_e(R)}$. After identifying an ordered basis of $V^{\oplus a_e(R)} \oplus V^{\oplus \rank_R(F^e_*R)-a_e(R)}$, the matrix of linear maps
    \[
    T_{F^e_*x,V}: V^{\oplus a_e(R)} \oplus V^{\oplus \rank_R(F^e_*R)-a_e(R)} \xrightarrow{\begin{pmatrix}A_{1,1,V} & A_{1,2,V} \\ {A_{2,1,V}} & A_{2,2,V}\end{pmatrix}} V^{\oplus a_e(R)} \oplus V^{\oplus \rank_R(F^e_*R)-a_e(R)}
    \]
    is so that all entries of the first $a_e(R)$ rows are in $IV$. The standard determinant expansion formula implies $\det(T_{F^e_*x,V}) \in I^{a_e(R)}V$
    and
    \[
    \nu(\det(T_{F^e_*x,V}))\geq a_e(R)\nu(I).
    \]
    But, $\left(T_{F^e_*x,V}\right)^{p^e} = x\id_{V^{\oplus \rank_R(F^e_*R)}}$ implies
    \[
    \det\left(T_{F^e_*x,V}\right)^{p^e} = \det\left(\left(T_{F^e_*x,V})^{p^e}\right)\right) = \det\left(x\id_{V^{\oplus \rank_R(F^e_*R)}}\right) = x^{\rank_R(F^e_*R)}.
    \]
    Therefore
    \[
    \rank_R(F^e_*R)\nu(x) = p^e\nu\left(\det\left(T_{F^e_*x,V})\right)\right)\geq p^ea_e(R)\nu(I).
    \]
    Dividing by $\rank_R(F^e_*R)$ provides the claimed inequality when $e_1=0$,
    \[
    \nu(x)\geq \frac{p^ea_e(R)}{\rank_R(F^e_*R)}\nu(I).
    \]

    When $e_1\geq 1$, $x\in I_e(I)$ implies $x^{p^{e_1}}\in I_{e+e_1}(I)$. By the above,
    \[
    p^{e_1}\nu(x) = \nu(x^{p^{e_1}}) \geq \frac{p^{e+e_1}a_{e+e_1}(R)}{\rank_R(F^{e+e_1}_*R)}\nu(I).
    \]
    Dividing by $p^{e_1}$ is the claimed inequality.
\end{proof}

\begin{corollary}
    \label{cor: valuation bound general}
    Let $(R,\fm,k)$ be either an $F$-finite local strongly $F$-regular domain of prime characteristic $p>0$ or a local normal domain essentially of finite type over $\CC$ of klt type. Let $\nu$ be a divisorial valuation of $R$. Then for every ideal $0\not=I\subseteq R$, $n\in\NN$, and $x\in I^{(n)}$, if $h$ is the maximal height of an associated prime of $I$,
    \[
    \nu(x)\geq \frac{ns(R)}{h}\nu(I).
    \]
\end{corollary}

\begin{proof}
    If $\fm$ is an associated prime of $I$, then $I^{(n)} = I^n$ and for every $x\in I^{n}$ it is clear that $\nu(x)\geq n\nu(I)\geq \frac{ns(R)}{h}\nu(I)$. In what follows, we assume that $\fm$ is not associated to $I$ and let $x\in I^{(n)}$. 
    
    First consider the case $(R,\fm,k)$ is an $F$-finite local strongly $F$-regular domain of prime characteristic $p>0$. For each $e\in\NN$ write $hp^e = q_en+r_e$ with $0\leq r_e<n$. Fix a nonzero element $0\not = g\in I^{(n)}$. Then $gx^{q_e}\in I^{(hp^e)}$. By \cref{lem: symbolic powers contained in splitting ideals}, $gx^{q_e}\in I_e(I)$. By \cref{thm: splitting ideals and valuation rings},
    \[
    \nu(g) + q_e\nu(x)\geq p^es(R)\nu(I).
    \]
    Multiplying by $n$,
    \[
    \nu(g)n + hp^e\nu(x)\geq \nu(g)n + q_en\nu(x) \geq p^ens(R)\nu(I).
    \]
    Dividing by $p^e$ and taking a limit as $e\to \infty$ gives $h\nu(x)\geq s(R)n\nu(I)$ and established the claimed inequality in the prime characteristic case.

    Now consider the case $(R,\fm,k)$ is a local normal domain essentially of finite type over $\CC$ and of klt type. Let $J\subseteq R$ be a choice of ideal containing $\nu$ as one of its Rees valuations. Choose a simultaneous log-resolution $\pi:X\to\Spec(R)$ of $I$ and $J$ that factors through the normalized blowups of $I$ and $J$. Let $E$ be a prime divisor of $X$ corresponding to $\nu$.

    Choose a suitable model of $R$, $I$, $x$, $\pi$, and $E$ over a finitely generated $\ZZ$-subalgebra $B\subseteq\CC$. Following the reduction to prime characteristic argument of \cite[proof of Theorem~4.4]{HochsterHuneke2002}, after shrinking $\Spec(B)$ as necessary, we may assume that for every general closed point $\mu\in\Spec(B)$, the element $x_\mu$ belongs to $I_\mu^{(n)}$ and the maximal height $h_\mu$ of an associated prime of $I_\mu$ is at most $h$. Moreover, we may assume that the reduction $E_\mu$ determines a divisorial valuation $\nu_\mu$ of $R_\mu$ so that $\nu_\mu(x_\mu)=\nu(x)$ and $\nu_\mu(I_\mu)=\nu(I)$. Indeed, these equalities follow by preserving the orders of vanishing of $x$ and $I$ along $E$ under reduction to prime characteristic.

    By \cref{cor: positive lim F-signature}, for a fixed $0<\varepsilon\ll 1$ we may further assume $s(R_\mu)>s(R)-\varepsilon$ for every general closed point $\mu\in\Spec(B)$. By the prime characteristic case, we obtain
    \[
        \nu(x)=\nu_\mu(x_\mu)
        \geq \frac{s(R_\mu)}{h_\mu}n\nu_\mu(I_\mu)
        \geq \frac{s(R)-\varepsilon}{h}n\nu(I).
    \]
    Allowing $\varepsilon\to 0$ shows that $\nu(x)\geq\frac{s(R)}{h}n\nu(I)$. 
\end{proof}

\begin{theorem}[Uniform Symbolic Topology Property]
    \label{cor: SFR rings enjoy USTP}
    Let $(R,\fm,k)$ be either an $F$-finite local strongly $F$-regular domain of prime characteristic $p>0$ or a local normal domain essentially of finite type over $\CC$ of klt type. Let $R$ have Krull dimension $d\geq 1$. Then for every ideal $0\not=I\subseteq R$, if $h$ is the maximal height of an associated prime of $I$, then for every $n\in\NN$
    \[
    I^{(n)}\subseteq \overline{I^{\left\lfloor \frac{s(R)}{h}n\right\rfloor}}, \quad\quad I^{\left(\left\lceil\frac{h}{s(R)}\right\rceil n\right)}\subseteq \overline{I^n}, \quad\mbox{and}\quad I^{\left(\left\lceil\frac{h}{s(R)}\right\rceil dn\right)}\subseteq I^n.
    \]
\end{theorem}

\begin{proof}
    Let $0\not=x\in I^{(n)}$. By \cref{cor: valuation bound general}, for each discrete valuation $\nu$ of $R$,
    \[
    \nu(x)\geq \frac{s(R)}{h}n\nu(I)\geq \left\lfloor\frac{s(R)}{h}n\right\rfloor \nu(I).
    \]
    By Rees valuative criteria for integral closure, $x\in \overline{I^{\left\lfloor \frac{s(R)}{h}n\right\rfloor }}$.

    If $x\in I^{\left(\left\lceil\frac{h}{s(R)}\right\rceil n\right)}$ then
    \[
    x\in \overline{I^{\left\lfloor\frac{s(R)}{h}\left\lceil\frac{h}{s(R)}n\right\rceil \right\rfloor}}\subseteq \overline{I^{\left\lfloor n\right\rfloor}} = \overline{I^n}.
    \]
    Similarly, if $x\in I^{\left(\left\lceil\frac{h}{s(R)}\right\rceil dn\right)}$ then $x\in \overline{I^{dn}}$. By \cite[Theorem 13.2.1 and Corollary 13.2.2]{HunekeSwanson2006}, $x\in I^n$.
\end{proof}

Let $k$ be a perfect field of characteristic $p>0$, and let $R_{d,n}$ denote the localization of the $n$th Veronese subring of $k[x_1,\ldots,x_d]$ with $d\geq 2$ at its homogeneous maximal ideal $\fm$. It is well known that $e(R_{d,n})=n^{d-1}$, while $s(R_{d,n})=1/n$ by \cite[Example 8 Page 320]{Singh2005}. In particular, $e(R_{d,n})=s(R_{d,n})^{-(d-1)}$. Consequently, for $d\geq 2$ and $\lambda<d-1$, there is no constant $C_d$ depending only on $d$ such that $e(R)\leq C_d/s(R)^\lambda$ for every $d$-dimensional $F$-finite strongly $F$-regular local domain $R$. Thus, the following corollary establishes an upper bound for multiplicity with the optimal exponent of the $F$-signature.

\begin{theorem}
    \label{cor: multiplicity bound}
    Let $(R,\fm,k)$ be either a local $F$-finite strongly $F$-regular domain of prime characteristic $p>0$ or an essentially finite type $\CC$-algebra of klt type. Let $R$ have Krull dimension $d\geq 1$. Then
    \[
    e(R)\leq \frac{d!}{s(R)^{d-1}}.
    \]
\end{theorem}

\begin{proof}
    First assume $R$ is of prime characteristic. Let $e\gg 0$ so that $\lfloor s(R)p^e\rfloor>0$. By \cref{cor: SFR rings enjoy USTP}, for every $e\in\NN$, $I_e(\fm)\subseteq \overline{\fm^{\lfloor s(R)p^e\rfloor}}.$
    Therefore
    \begin{align*}
    \frac{a_e(R)}{\rank_R(F^e_*R)} &= \frac{\mathcal{L}(R/I_e(\fm))}{p^{ed}}\\
    &\geq \frac{\mathcal{L}\left(R/\overline{\fm^{\lfloor s(R)p^e\rfloor}}\right)}{p^{ed}}\\
    &=\frac{d!\mathcal{L}\left(R/\overline{\fm^{\lfloor s(R)p^e\rfloor}}\right)}{\lfloor s(R)p^e\rfloor^d}\times \frac{\lfloor s(R)p^e\rfloor^d}{d!p^{ed}}.
    \end{align*}
    Consequently,
    \begin{align*}
        s(R) &= \lim_{e\to \infty} \frac{a_e(R)}{\rank_R(F^e_*R)}\\
        &\geq \lim_{e\to \infty} \left(\frac{d!\mathcal{L}\left(R/\overline{\fm^{\lfloor s(R)p^e\rfloor}}\right)}{\lfloor s(R)p^e\rfloor^d}\times \frac{\lfloor s(R)p^e\rfloor^d}{d!p^{ed}}\right)\\
        &= \frac{e(R)s(R)^d}{d!}.
    \end{align*}
    Rearranging terms,
    \[
    e(R)\leq \frac{d!}{s(R)^{d-1}}.
    \]

    The characteristic $0$ inequality follows by standard reduction to positive characteristic arguments outlined in previous arguments.
\end{proof}

\begin{corollary}
    \label{cor: embedding dimension bound}
    Let $(R,\fm,k)$ be either a local $F$-finite strongly $F$-regular domain of prime characteristic $p>0$ or an essentially finite type $\CC$-algebra of klt type. Let $R$ have Krull dimension $d\geq 1$. Then
    \[
    \edim(R) \leq d-1 + \frac{d!}{s(R)^{d-1}}.
    \]
\end{corollary}
\begin{proof}
    Abhyankar's inequality is that $\edim(R) - d + 1\leq e(R)$. By \cref{cor: multiplicity bound},
    \[
    \edim(R)\leq d-1 + e(R) \leq d-1 + \frac{d!}{s(R)^{d-1}}.\qedhere
    \]
\end{proof}

By \cite[Main Theorem~2]{PolstraUSTP}, if $(R,\fm,k)$ is a local $F$-finite strongly $F$-regular domain of prime characteristic $p>0$, then there exists $e_0\in\NN$ so that for every ideal $I\subseteq R$ and $e\geq 1$, $I_{e+e_0}(I)\subseteq I^{[p^e]}$. The methods used in \cite{PolstraUSTP} only provide the existence of such a constant. Here we able to provide an explicit choice of $e_0\in\NN$ dependent only upon the dimension and $F$-signature of $R$. 

\begin{theorem}
    \label{thm: splitting ideal containments in SFR rings}
    Let $(R,\fm,k)$ be a local $F$-finite strongly $F$-regular domain of prime characteristic $p>0$ and Krull dimension $d\geq 1$. Let $0\not=I\subseteq R$ be an ideal. For every $e\geq 1$,
    \[
    I_e(I)\subseteq \overline{I^{\left\lfloor p^es(R)\right\rfloor}}.
    \]
    Moreover, if $e_0 = \left\lceil \log_p\left(\frac{d}{s(R)}\right)\right\rceil$, then for every $e\geq 1$,
    \[
    I_{e+e_0}(I)\subseteq I^{[p^e]}.
    \]
\end{theorem}

\begin{proof}
    let $0\not=x\in I_e(I)$. By \cref{thm: splitting ideals and valuation rings}, for each discrete valuation $\nu$ of $R$, $\nu(x)\geq p^es(R)\nu(I)\geq \lfloor p^es(R)\rfloor\nu(I)$.  By Rees valuative criteria for integral closure, \cite[Theorem 10.2.2(3)]{HunekeSwanson2006}, $x\in \overline{I^{\left\lfloor p^es(R)\right\rfloor}}$.
    
    Let $e_0 = \left\lceil \log_p\left(\frac{d}{s(R)}\right)\right\rceil$. Then for every $e\in\NN$, $I_{e+e_0}(I)\subseteq \overline{I^{\lfloor p^{e+e_0}s(R)\rfloor}}.$ By our choice of $e_0$, $p^{e_0}s(R)\geq \frac{d}{s(R)}s(R)=d$. Therefore $I_{e+e_0}(I)\subseteq \overline{I^{dp^e}}\subseteq I^{[p^e]}$ by \cite[Theorem 13.2.1 and Corollary 13.2.2]{HunekeSwanson2006}.
\end{proof}

\begin{corollary}
    \label{cor: F-sig lower bound for alpha F invariant}
    Let $(R,\fm,k)$ be a local $F$-finite and $F$-pure domain of prime characteristic $p>0$. Then $s(R)\leq \alpha_F(R)$.
\end{corollary}

\begin{proof}
    The content of the corollary is trivial if $R$ is not strongly $F$-regular. If $R$ is strongly $F$-regular, by \cref{thm: splitting ideal containments in SFR rings}, for each $e\geq 1$, $I_e(\fm)\subseteq \overline{\fm^{\lfloor s(R)p^e\rfloor}}$. Equivalently, $m_e(R)\geq \lfloor s(R)p^e\rfloor$. Dividing by $p^e$ and taking a limit as $e\to\infty$ gives
    \[
    \alpha_F(R) = \lim_{e\to\infty}\frac{m_e(R)}{p^e} \geq \lim_{e\to \infty}\frac{\lfloor s(R)p^e\rfloor}{p^e} = s(R).\qedhere
    \]
\end{proof}

In characteristic zero, the notion of normalized volume is motivated by \cite[Theorem~1.2]{Li18}. Li established an Izumi-type inequality whose dependence on the valuation is controlled by log discrepancy. Specifically, if $(R,\fm,k
)$ is the local ring of a complex klt singularity, then there exists a constant $C=C(R)>0$ such that
\[
    \nu(\fm)\ord_{\fm}(x)
    \leq \nu(x)
    \leq C A(\nu)\ord_{\fm}(x)
\]
for every real valuation $\nu$ centered on $\fm$ with finite log discrepancy and every $0\not=x\in R$. Applying these inequalities to the powers of $x$ and passing to the limit gives the same inequalities with $\ord_{\fm}(x)$ replaced by $\overline{\ord}_{\fm}(x)$. In \cref{cor: discrepancy bound}, we simply observe that a prime characteristic analogue of Li's valuation estimate theorem is true in $\QQ$-Gorenstein strongly $F$-regular singularities with the singularity-dependent constant given explicitly by either the reciprocal of the $\alpha_F$-invariant or the $F$-signature of $R$.

\begin{corollary}
    \label{cor: discrepancy bound}
    Let $(R,\fm,k)$ be an $F$-finite strongly $F$-regular local $\QQ$-Gorenstein domain of prime characteristic $p>0$. Then for every $0\not=x \in \fm$ and divisorial valuation $\nu$ centered on $\fm$, if $A(\nu)$ is the log-discrepancy of $\nu$, then
    \[
    \nu(\fm)\overline{\ord}_{\fm}(x) \leq \nu(x)\leq \frac{A(\nu)}{\alpha_F(R)}\overline{\ord}_{\fm}(x)  \leq \frac{A(\nu)}{s(R)}\overline{\ord}_{\fm}(x).
    \]
\end{corollary}
\begin{proof}
    It is clear that $\nu(\fm)\overline{\ord}_{\fm}(x)\leq \nu(x)$. Let $t<\FPT_R(x)$ be a rational number so that the pair $(R,x^t)$ is strongly $F$-regular. By \cite[Theorem~3.3]{HaraWatanabe2002}, $A(\nu)>t\nu(x)$. Therefore $A(\nu)\geq \fpt_R(x)\nu(x)\geq \frac{\alpha_F(R)}{\overline{\ord}_{\fm}(x)}\nu(x)$. Solving for $\nu(x)$ and the inequality of \cref{cor: F-sig lower bound for alpha F invariant} provides the claimed inequalities.
\end{proof}

\section{Izumi-Rees Bounds}

Let $R$ be an excellent domain of finite Krull dimension $d$.
Then $R$ is known to enjoy the Uniform Brian\c{c}on-Skoda Property; there exists a constant $b\in\NN$ such that for every ideal $I\subseteq R$ and every $n\in\NN$,
\[
    \overline{I^{n+b}}\subseteq I^n.
\]
See \cite{MaMcDonaldRGSchwede2026,Huneke1992}.

Effective uniform Brian\c{c}on-Skoda bounds are known under additional hypotheses. Classically, if $(R,\fm,k)$ is a $d$-dimensional ring with pseudo-rational singularities, e.g. $R$ is prime characteristic, excellent, and $F$-rational by \cite{Smith1997}, then $b=d-1$ is a Uniform Brian\c{c}on-Skoda bound by \cite{LipmanSathaye1981}. If $(R,\fm,k)$ is a $d$-dimensional $F$-pure local ring of prime characteristic $p>0$, then $b=d$ is a uniform Brian\c{c}on-Skoda bound by \cite[Corollary~3.3]{WheelerZhang2025}. The same bound holds for Cohen-Macaulay $F$-injective local rings \cite[Corollary~3.12]{WheelerZhang2025}. More generally, Ma, McDonald, R.G., and Schwede prove that $b=d$ is a uniform Brian\c{c}on-Skoda bound for Du Bois singularities and $F$-pure singularities \cite[Theorem~3.18]{MaMcDonaldRGSchwede2026}.

The results of this section uses effective Uniform Brian\c{c}on-Skoda bounds to derive an effective Uniform Artin-Rees bound of a principal ideal $(x)\subseteq R$ as a function of the $F$-pure threshold of $x$.

\begin{lemma}
    \label{lem: valution bound via fpt}
    Let $(R,\fm,k)$ be an $F$-finite and $F$-pure domain of prime characteristic $p>0$ and Krull dimension $d$. Let $0\not=x\in \fm$ have the property $\fpt(x)>0$. Then for every ideal $I\subseteq R$ and $0\not=y\in R$,
    \[
    \overline{\ord}_I(xy)\leq \ord_I(y) + \frac{d+1}{\fpt(x)}.
    \]
    If $R$ is strongly $F$-regular, then
    for every $y\in R$,
    \[
    \overline{\ord}_I(xy)\leq \ord_I(y) + \frac{d}{\fpt(x)}.
    \]
\end{lemma}
\begin{proof}
    Let
    \[
    b = \begin{cases}
        d & \mbox{if $R$ is strongly $F$-regular}\\
        d+1 & \mbox{if $R$ is not strongly $F$-regular}
    \end{cases}.
    \]
    The constant $b$ has the property that for ideals $I\subseteq R$ and $n\in\NN$, $\overline{I^{n+b-1}}\subseteq I^n$. 
    
    Let $\nu_e(x) = \max\{n\in\NN\mid x^n\not\in I_e(\fm)\}$ so that $\FPT(x) = \lim_{e\to\infty}\frac{\nu_e(x)}{p^e}$. We first claim $x^{\nu_e(x)}y^{p^e}\not\in \overline{I^{p^e(\ord_{I}(y)+b)}}$. Otherwise,
    \begin{align*}
        x^{\nu_e(x)}y^{p^e} &\in \overline{I^{p^e(\ord_{I}(y)+b)}} \\ 
        &= \overline{\left(I^{[p^e]}\right)^{\ord_I(y)+b}}\\
        &\subseteq \left(I^{[p^e]}\right)^{\ord_I(y)+1}\\
        &= \left(I^{\ord_I(y)+1}\right)^{[p^e]}.
    \end{align*}
    If $\varphi\in\Hom_R(F^e_*R,R)$ is so that $\varphi(F^e_*x^{\nu_e(x)})=1$ then $\varphi(F^{e}_*x^{\nu_e(x)}y^{p^e}) = y\in I^{\ord_{I}(y)+1}$, a contradiction. 
    
    Since $x^{\nu_e(x)}y^{p^e} = (xy)^{\nu_e(x)}y^{p^e-\nu_e(x)}\not\in \overline{I^{p^e(\ord_{I}(y)+b)}}$ we necessarily have 
    \[
    \nu_e(x)\overline{\ord}_{I}(xy) + (p^e-\nu_e(x))\overline{\ord}_{I}(y)<p^e(\ord_{I}(y)+b).
    \]
    Consequently,
    \[
    \nu_e(x)\overline{\ord}_{I}(xy)<\nu_e(x)\ord_{I}(y) + p^eb.
    \]
    Dividing by $p^e$ and taking a limit as $e\to \infty$,
    \[
    \fpt(x)\overline{\ord}_{I}(xy)\leq \fpt(x)\ord_{I}(y) + b,
    \]
    i.e., $\overline{\ord}_{I}(xy)\leq \ord_{I}(y) +\frac{b}{\FPT(x)}$.
\end{proof}

\begin{corollary}
    \label{cor: Artin Rees bound via Fpure threshold}
    Let $(R,\fm,k)$ be an $F$-finite and $F$-pure domain of prime characteristic $p>0$ and Krull dimension $d$. Let $0\not=x\in \fm$ satisfy $\fpt(x)>0$. Then for every ideal $I\subseteq R$ and $n\in\NN$,
    \[
    \left(\overline{I^{n}}:_Rx\right)\subseteq I^{n-\left\lfloor \frac{d+1}{\FPT(x)}\right\rfloor}.
    \]
    If $R$ is strongly $F$-regular, then
    for every ideal $I\subseteq R$ and $n\in\NN$,
    \[
    \left(\overline{I^{n}}:_Rx\right)\subseteq I^{n-\left\lfloor \frac{d}{\FPT(x)}\right\rfloor}.
    \]
\end{corollary}

\begin{proof}
    Set
    \[
        b = \begin{cases}
        d & \mbox{if $R$ is strongly $F$-regular}\\
        d+1 & \mbox{if $R$ is not strongly $F$-regular}
        \end{cases}.
    \]
    
    Let $0\not=y\in R$. By \cref{lem: valution bound via fpt}, $\overline{\ord}_{I}(xy)\leq \ord_{I}(y) +\frac{b}{\FPT(x)}$. Since $\ord_I(y)$ is an integer,
    \[
        \ord_I(y)\geq n - \left\lfloor \frac{b}{\fpt(x)} \right\rfloor.
    \]
    and $y\in I^{n-\left\lfloor \frac{b}{\fpt(x)}\right\rfloor}$.
\end{proof}

\begin{lemma}\label{lem: Artin Rees bound via log canonical threshold}
    Let $(R,\fm,k)$ be a $d$-dimensional normal domain essentially of finite type over $\bC$. Let $X=\Spec(R)$ and $\Delta\geq 0$ be an effective $\QQ$-divisor so that $K_X+\Delta$ is $\QQ$-Cartier and $(X,\Delta)$ is klt. Then for every proper nonzero ideal $I \subseteq R$, $0\not=x\in\fm$, and $0\not=y\in R$,
    \[
    \overline{\ord}_I(xy)\leq \ord_I(y) + \frac{d}{\lct_{(X,\Delta)}(x)}.
    \]
\end{lemma}
\begin{proof}
    Fix $0\neq y\in R$. Suppose, for contradiction, that
    $\overline{\ord}_I(xy)-\ord_I(y)>\frac{\ell(I)}{\lct_{(X,\Delta)}(x)}$,
    where $\ell(I)$ is the analytic spread of $I$.
    Thus, we may choose positive integers $a<b$ such that
    $\frac{a}{b}<\lct_{(X,\Delta)}(x)$ and
    $a(\overline{\ord}_I(xy)-\ord_I(y))>b\ell(I)$.
    Consequently,
    \begin{align*}
        \overline{\ord}_I(x^ay^b)
        &=\overline{\ord}_I((xy)^ay^{b-a})\\
        &\geq a\overline{\ord}_I(xy)+(b-a)\ord_I(y)\\
        &=b\ord_I(y)+a(\overline{\ord}_I(xy)-\ord_I(y))\\
        &>b(\ord_I(y)+\ell(I)).
    \end{align*}
    In particular,
    $x^ay^b\in\overline{I^{b(\ord_I(y)+\ell(I))}}$.

    Let $v$ be any divisorial valuation over $X=\Spec R$.
    By the valuative criterion for integral closure,
    $\frac{a}{b}v(x)+v(y)\geq(\ord_I(y)+\ell(I))v(I)$.
    Since $\frac{a}{b}<\lct_{(X,\Delta)}(x)$, the pair
    $(X,\Delta + \frac{a}{b}\Div(x))$ is klt.
    Thus, $A_{(X,\Delta)}(v)>\frac{a}{b}v(x)$, and hence
    \[
        A_{(X,\Delta)}(v)+v(y)>(\ord_I(y)+\ell(I))v(I).
    \]
    Since this holds for every divisorial valuation over $X$,
    we conclude that
    $y\in\mathcal{J}((X,\Delta),I^{\ord_I(y)+\ell(I)})$.
    By \cite[Corollary 1.4]{LLS2008},
    $\mathcal{J}((X,\Delta),I^{\ord_I(y)+\ell(I)})
    \subseteq I^{\ord_I(y)+1}$.
    Thus, $y\in I^{\ord_I(y)+1}$, a contradiction.

    Therefore,
    \[
        \overline{\ord}_I(xy)-\ord_I(y)
        \leq\frac{\ell(I)}{\lct_{(X,\Delta)}(x)}
        \leq\frac{d}{\lct_{(X,\Delta)}(x)}.\qedhere
    \]
\end{proof}

\begin{corollary}\label{cor: Artin Rees bound via log canonical threshold}
    Let $(R,\fm,k)$ be a $d$-dimensional normal domain essentially of finite type over $\bC$. Let $X=\Spec(R)$ and $\Delta\geq 0$ be an effective $\QQ$-divisor so that $K_X+\Delta$ is $\QQ$-Cartier and $(X,\Delta)$ is klt. Let $0 \neq x \in \fm$. Then for every proper nonzero ideal $I \subseteq R$ and every $n \geq 1$,
    \[
    \left(\overline{I^{n}}:_Rx\right)\subseteq \overline{I^{n-\left\lfloor \frac{d}{\lct_{(X,\Delta)}(x)}\right\rfloor}}.
    \]
\end{corollary}

\begin{proof}
    This is an application of \cref{lem: Artin Rees bound via log canonical threshold} and follows by the argument presented in \cref{cor: Artin Rees bound via Fpure threshold}.
\end{proof}

The extension of the preceding result to log canonical singularities analogous to $F$-pure singularities described in \cref{lem: valution bound via fpt} and \cref{cor: Artin Rees bound via Fpure threshold} presents two difficulties. First, log canonical pairs may admit divisorial valuations with vanishing log discrepancy, and consequently the log canonical threshold of a nonzero element may vanish. Second, the proof relies on strict positivity of log discrepancies to establish membership in a multiplier ideal, a condition that need not hold for log canonical pairs. Although log canonical singularities satisfy a Brian\c{c}on–Skoda theorem with an additional shift of one, analogous to the distinction between the strongly $F$-regular and $F$-pure cases, this does not directly yield an estimate in terms of log canonical thresholds. Thus, obtaining a meaningful extension requires either additional hypotheses on the log canonical centers or a different approach to measuring the singularities.

\begin{theorem}[Uniform Izumi-Rees Property]
    \label{thm: Uniform Izumi-Rees for SFR}
    Let $(R,\fm,k)$ be either a local $F$-finite strongly $F$-regular domain of prime characteristic $p>0$ or an essentially of finite type $\CC$-algebra of klt type. Let $R$ have Krull dimension $d\geq 1$. Let $0\not=x\in \fm$ and $I\subseteq R$ an $\fm$-primary ideal.
    \begin{enumerate}
        \item If $R$ has prime characteristic,
        \[
        e_I(R/xR)\leq \left\lfloor \frac{d}{\FPT(x)}\right\rfloor e_I(R)\leq \frac{de_I(R)}{\alpha_F(R)}\overline{\ord}_{\fm}(x) \leq \frac{de_I(R)}{s(R)}\overline{\ord}_{\fm}(x).
        \]
        \item If $R$ has characteristic $0$ and $\Delta\geq 0$ an effective boundary $\QQ$-divisor so that $(R,\Delta)$ is klt,
        \[
        e_I(R/xR)\leq \left\lfloor \frac{d}{\lct_{(R,\Delta)}(x)}\right\rfloor e_I(R)\leq  \frac{de_I(R)}{s(R)}\overline{\ord}_{\fm}(x).
        \]
    \end{enumerate}
    
\end{theorem}

\begin{proof}
    First consider the case $R$ has prime characteristic. For each $n\in\NN$ there is a short exact sequence
    \[
    0\to \frac{R}{(I^{n}:_Rx)}\xrightarrow{\cdot x}\frac{R}{I^n}\to \frac{R}{(I^n,x)}\to 0.
    \]
    By \cref{cor: Artin Rees bound via Fpure threshold}, for each $n\in\NN$, $(I^n:_Rx)\subseteq (\overline{I^n}:_Rx)\subseteq I^{n-\left\lfloor \frac{d}{\FPT(x)}\right\rfloor}$. Therefore
    \begin{align*}
        \mathcal{L}\left(\frac{R}{(I^n,x)}\right) &= \mathcal{L}\left(\frac{R}{I^n}\right) -\mathcal{L}\left(\frac{R}{(I^n:x)}\right)\\
        &\leq \mathcal{L}\left(\frac{R}{I^n}\right) -\mathcal{L}\left(\frac{R}{I^{n-\left\lfloor \frac{d}{\FPT(x)}\right\rfloor}}\right).
    \end{align*}
    The functions on the left and right hand side are polynomial expressions in $n$ of degree $d-1$ for $n\gg 0$. The leading coefficient of the left hand side polynomial function is $\frac{e_I(R/xR)}{(d-1)!}$ and the leading coefficient of the right hand side polynomial function is $\frac{e_I(R)}{(d-1)!}\left\lfloor \frac{d}{\FPT(x)}\right\rfloor$. Therefore
    \[
    e_I(R/xR)\leq e_I(R)\left\lfloor \frac{d}{\FPT_R(x)}\right\rfloor\leq \frac{de_I(R)}{\FPT_R(x)}.
    \]
    Recall that $\alpha_F(R) = \inf\{\FPT_R(y)\overline{\ord}_{\fm}(y)\mid 0\not=y\in \fm\}$, implying $\FPT_R(x)\geq \frac{\alpha_F(R)}{\overline{\ord}_{\fm}(x)}$. By \cref{cor: F-sig lower bound for alpha F invariant}, $\fpt_R(x)\geq \frac{\alpha_F(R)}{\overline{\ord}_{\fm}(x)}\geq \frac{s(R)}{\overline{\ord}_{\fm}(x)}$. Therefore
    \[
    e_I(R/xR)\leq \left\lfloor \frac{d}{\FPT(x)}\right\rfloor e_I(R)\leq  \frac{de_I(R)}{\alpha_F(R)}\overline{\ord}_{\fm}(x) \leq \frac{de_I(R)}{s(R)}\overline{\ord}_{\fm}(x).
    \]

    Now suppose $R$ has characteristic $0$. Let $0<t<\lct_{(R,\Delta)}(x)$ be rational and $0<\varepsilon<s(R)$. Since $(R,\Delta+t\operatorname{div}(x))$ is klt, by \cite[Corollary~3.4]{Takagi2004}, we may choose a suitable model over a finitely generated $\ZZ$-subalgebra $B\subseteq\CC$ so that $(R_\mu,\Delta_\mu+t\operatorname{div}(x_\mu))$ is strongly $F$-regular for every general closed point $\mu\in\Spec(B)$. In particular, $\FPT_{R_\mu}(x_\mu)\geq t$. By the definition of characteristic zero $F$-signature, we may further assume $s(R_\mu)>s(R)-\varepsilon$. After enlarging the model and shrinking $\Spec(B)$ as necessary, we may also assume $\dim(R_\mu)=d$, $e_{I_\mu}(R_\mu)=e_I(R)$, and $e_{I_\mu}(R_\mu/x_\mu R_\mu)=e_I(R/xR)$.

    Choosing a log resolution of $\fm$ and spreading out this resolution together with the preceding data, we may further assume $\overline{\ord}_{\fm_{\mu}}(x) = \overline{\ord}_{\fm}(x)$ for a general closed point $\mu$ of the coefficient ring of the model. Applying the prime characteristic case to $R_\mu$, letting $t$ approach $\lct_{(R,\Delta)}(x)$ and $\varepsilon$ approach $s(R)$ provides the desired inequalities.
\end{proof}

\begin{corollary}[Izumi-Rees Bounds]
    \label{cor: explicit Izumi-Rees bound}
     Let $(R,\fm,k)$ be either a local $F$-finite strongly $F$-regular domain of prime characteristic $p>0$ or an essentially of finite type $\CC$-algebra of klt type. Let $R$ have Krull dimension $d\geq 1$. Let $I\subseteq R$ be an $\fm$-primary ideal, $\nu_1,\nu_2$ Rees valuations and of $I$, and $0\not=x\in R$.
     
     \begin{enumerate}
         \item If $R$ has prime characteristic,
         \[
         \nu_1(x)\leq \frac{de_I(R)}{\alpha_F(R)} \nu_2(x) \leq \frac{de_I(R)}{s(R)}\nu_2(x).
         \]
         \item If $R$ has characteristic $0$,
         \[
         \nu_1(x)\leq \frac{de_I(R)}{s(R)}\nu_2(x).
         \]
     \end{enumerate}
\end{corollary}

\begin{proof}
    We may assume $x\in \fm$, else the claimed inequalities are trivial. By a theorem of Rees, \cite[Main Result]{Rees1961} see also \cite[Theorem~2.8]{Polstra2026} and \cite[Corollary 4.13]{CutkoskyMontano2025}, if $\{\nu_1,\ldots,\nu_t\}$ are the divisorial valuations associated to the normalized blowup of $I$, then there are positive integers $d_1,\ldots,d_t$ so that for every $0\not=y\in\fm$,
    \[
    e_I(R/yR) = \sum_{i=1}^td_i\nu_i(y).
    \]
    Moreover, by \cite[Corollary~2.9]{Polstra2026},
    \[
    e_I(R) = \sum_{i=1}^t d_i\nu_i(I).
    \]

    Let $0\not=x\in R$. The claimed inequality is trivial if $x\not\in\fm$. If $0\not=x\in\fm$, then
    \[
    \nu_1(x)\leq \sum_{i=1}^td_i\nu_i(x) = e_I(R/xR).
    \]
    By \cref{thm: Uniform Izumi-Rees for SFR}, if $R$ is of prime characteristic,
    \[
    \nu_1(x)\leq \frac{de_I(R)}{\alpha_F(R)}\overline{\ord}_{\fm}(x) \leq \frac{de_I(R)}{s(R)}\overline{\ord}_{\fm}(x),
    \]
    and if $R$ has characteristic $0$,
    \[
    \nu_1(x)\leq \frac{de_I(R)}{s(R)}\overline{\ord}_{\fm}(x).
    \]
    Because $x\in\overline{\fm^{\overline{\ord}_{\fm}(x)}}$, $\overline{\ord}_{\fm}(x)\nu_2(\fm)\leq \nu_2(x)$. Therefore in prime characteristic
    \[
    \nu_1(x)\leq \frac{de_I(R)}{\alpha_F(R)}\frac{\nu_2(x)}{\nu_2(\fm)}\leq \frac{de_I(R)}{\alpha_F(R)}\nu_2(x)\leq \frac{de_I(R)}{s(R)}\nu_2(x),
    \]
    and in characteristic $0$,
    \[
    \nu_1(x)\leq \frac{de_I(R)}{s(R)}\frac{\nu_2(x)}{\nu_2(\fm)}\leq \frac{de_I(R)}{s(R)}\nu_2(x)\qedhere
    \]
\end{proof}

\bibliographystyle{skalpha}
\bibliography{main}



\end{document}